\documentclass[11pt,reqno]{amsart}
\usepackage{fixltx2e}                    
\usepackage[utf8]{inputenc}            
\usepackage{amsmath}                     
\usepackage{amssymb, latexsym, stmaryrd, amsthm, dsfont, amsfonts, amsbsy,amsthm, amsmath, mathrsfs}            
\usepackage{mathtools}                   
\usepackage{bm}                          
\usepackage{enumerate}                   
\usepackage{verbatim}                    
\usepackage{url}                         
\usepackage[osf]{mathpazo}  

\newcommand\blfootnote[1]{%
  \begingroup
  \renewcommand\thefootnote{}\footnote{#1}%
  \addtocounter{footnote}{-1}%
  \endgroup
}

\usepackage{microtype}  
\usepackage[all, knot]{xy}
        \xyoption{arc} 
        \xyoption{web}                 
\makeatletter                            
\def\MT@register@subst@font{\MT@exp@one@n\MT@in@clist\font@name\MT@font@list
 \ifMT@inlist@\else\xdef\MT@font@list{\MT@font@list\font@name,}\fi}
\makeatother

\usepackage[pdftex,bookmarks,bookmarksnumbered,linktocpage,   
         colorlinks,linkcolor=blue,citecolor=blue]{hyperref}

\newcommand{\bit}{\begin{itemize}}    
\newcommand{\eit}{\end{itemize}}
\newcommand{\ben}{\begin{enumerate}}
\newcommand{\een}{\end{enumerate}}
\newcommand{\benormal}{\ben[\normalfont 1.]}   

\let\enormal\een
\newcommand{\benroman}{\ben[\normalfont (i)]}  

\newcommand{\bde}{\begin{description}}
\newcommand{\ede}{\end{description}}

\newcommand{\?}{\ensuremath{\mkern0.4\thinmuskip}}   

\let\leq=\leqslant
\let\geq=\geqslant
\let\Box=\square                            

\let\epsilon=\varepsilon
\let\Lambda\varLambda
\let\Gamma\varGamma
\let\Delta\varDelta
\let\Lambda\varLambda
\let\Omega\varOmega
\let\Theta\varTheta
\let\Xi\varXi
\let\Pi\varPi
\let\Sigma\varSigma

\let\clsys=\mathcal                             
\let\class=\mathsf                              
\let\oper=\mathbb                               

\bmdefine{\A}{A}                                
\bmdefine{\2}{2}
\bmdefine{\B}{B}
\bmdefine{\D}{D}
\bmdefine{\M}{M}                                
\bmdefine{\LLL}{L}                              
\bmdefine{\Fm}{Fm}                              
\bmdefine{\zerou}{[0{,}1]}  
\bmdefine{\T}{T}                                

\newcommand{\FF}{\clsys{F}}

\newcommand{\LL}{\mathcal{L}}                   

\newcommand{\VVV}{\oper{V}}                     

\newcommand{\Con}{\mathrm{Con}}                            
\newcommand{\Var}{\mathnormal{V\mkern-.8\thinmuskip ar}} 
\newcommand{\Th}{\clsys{T}\mkern-.5\thinmuskip\mathnormal{h}} 
\newcommand{\FFi}{\FF\mkern-.5\thinmuskip\mathnormal{i}}  
\newcommand{\FFiL}{\FFi_{\LL}}                            
\newcommand{\FFiLA}{\FFi_{\LL}\A}      
\newcommand{\Fi}{\text{\textsl{Fi}}}                      
\newcommand{\FiLA}{\Fi_{\LL}^{\A}}                        
\newcommand{\Mod}{\class{Mod}}
\newcommand{\ModS}{\class{Mod}^{\textup{Su}}}

\newcommand{\Alg}{\class{Alg}}

\newcommand{\LMod}{\class{LMod}}                 
\newcommand{\LModstar}{\LMod^{\boldstar}}                 
\newcommand{\LModstarL}{\LModstar\!\LL}                 
\newcommand{\ModSL}{\ModS\!\LL}

\bmdefine{\boldstar}{\mathchoice{\textstyle*}{\textstyle*}{\textstyle*}{\scriptstyle*}}
\newcommand{\Modstar}{\Mod^{\boldstar}}
\newcommand{\ModstarL}{\Modstar\!\LL}

\newcommand{\Algstar}{\Alg^{\boldstar}\!\?}
\newcommand{\AlgstarL}{\Algstar\!\?\LL}

\bmdefine{\btau}{\tau}                                  
\bmdefine{\brho}{\rho}                                  

\newcommand{\vdashL}{\vdash_{\!\LL}}                      
\newcommand{\sineq}{\mathrel{\dashv\mkern1.5mu\vdash}}  

\bmdefine{\leibniz}{\Omega}        

\bmdefine{\frege}{\Lambda}         

\makeatletter
\newcommand{\tarskidsp}{\mathord%
   {\m@th\raisebox{0pt}[0pt][0pt]{$\stackrel%
   {\raisebox{-2.7pt}[0ex][0pt]{$\displaystyle \,\?\thicksim$}}%
   {\displaystyle\leibniz}$}}}
\newcommand{\tarskitxt}{\mathord%
   {\m@th\raisebox{0pt}[0pt][0pt]{$\stackrel%
   {\raisebox{-2.7pt}[0ex][0pt]{$\,\?\thicksim$}}{\displaystyle\leibniz}$}}}
\newcommand{\tarskiscr}{\mathord%
   {{\m@th\raisebox{0pt}[0pt][0pt]{$\stackrel%
   {\raisebox{-2.4pt}[0ex][0pt]{$\scriptstyle \,\?\thicksim$}}%
   {\scriptstyle\leibniz}$}}}}
\newcommand{\tarskiscrscr}{\mathord%
   {{\m@th\raisebox{0pt}[0pt][0pt]{$\stackrel%
   {\raisebox{-2pt}[0ex][0pt]{$\scriptscriptstyle \,\?\thicksim$}}%
   {\scriptscriptstyle\leibniz}$}}}}
\newcommand{\tarski}{\@ifnextchar ^ %
   {\mathchoice{\tarskidsp\kern-.07em}{\tarskitxt\kern-.07em}%
   {\tarskiscr\kern-.07em}{\tarskiscrscr\kern-.07em}}%
   {\mathchoice{\tarskidsp}{\tarskitxt}{\tarskiscr}{\tarskiscrscr}}}
\makeatother

\theoremstyle{theorem}
\newtheorem{Theorem}{Theorem}[section]
\newtheorem{Lemma}[Theorem]{Lemma}
\newtheorem{Corollary}[Theorem]{Corollary}
\newtheorem{fact}[Theorem]{\textbf{Fact}}

\theoremstyle{definition}
\newtheorem{Definition}[Theorem]{Definition}
\newtheorem{exa}[Theorem]{Example}

\theoremstyle{remark}

\newtheorem{problem}{\bf Problem}
\newtheorem{Remark}[Theorem]{Remark}

\newcommand{\lo}{\mathcal{L}} 

\newcommand{\C}{\boldsymbol{C}} 

\allowdisplaybreaks[1]

\begin{document}
\title[A study of truth predicates in matrix semantics]{A study of truth predicates in matrix semantics}

\author{Tommaso Moraschini}
\email{moraschini@cs.cas.cz}
\address{Institute of Computer Science, Czech Academy of Sciences, Pod Vod\'{a}renskou v\v{e}\v{z}\'{i} 2, 182 07 Prague 8, Czech Republic}
\date{\today}
\maketitle

\begin{abstract}
Abstract algebraic logic is a theory that provides general tools for the algebraic study of arbitrary propositional logics. According to this theory, every logic $\LL$ is associated with a matrix semantics $\ModstarL$. This paper is a contribution to the systematic study of the so-called \textit{truth sets} of the matrices in $\ModstarL$. In particular, we show that the fact that the truth sets of $\ModstarL$ can be defined by means of equations with universally quantified parameters is captured by an order-theoretic property of the Leibniz operator restricted to deductive filters of $\LL$. This result was previously known for equational definability without parameters. Similarly, it was known that the truth sets of $\ModstarL$ are implicitly definable if and only if the Leibniz operator is injective on deductive filters of $\LL$ over every algebra. However, it was an open problem whether the injectivity of the Leibniz operator transfers from the theories of $\LL$ to its deductive filters over arbitrary algebras. We show that this is the case for logics expressed in a countable language, and that it need not be true in general. Finally we consider an intermediate condition on the truth sets in $\ModstarL$ that corresponds to the order-reflection of the Leibniz operator.
\end{abstract}

\section{Introduction}
\blfootnote{\textit{Keywords}: \textup{Abstract algebraic logic}, \textup{truth predicate}, \textup{equational definability}, \textup{truth-equational logic}, \textup{protoalgebraic logic}, \textup{Leibniz hierarchy}, \textup{Leibniz operator}, \textup{implicit definability}, \textup{matrix semantics}, \textup{algebraic semantics}, \textup{propositional logic}, \textup{protodisjunction}, \textup{protoconjunction}.} 
\blfootnote{\textup{2010} \textit{Mathematics Subject Classification}: \textup{03G27}, \textup{03G10}, \textup{03B22}.}

Abstract algebraic logic (AAL for short) is a theory that aims to provide general tools for the algebraic study of arbitrary propositional logics \cite{BP86,BP89,Cz01,AAL-AIT-f,FJa09,FJaP03b}. According to this theory, every (propositional) logic $\LL$ is associated with a matrix semantics $\ModstarL$ with respect to which $\LL$ is sound and complete. In the best-known cases the class of matrices $\ModstarL$ coincides with the intended algebraic semantics of $\LL$, e.g., in the case of superintuitionistic logics $\ModstarL$ is the class of matrices based on a variety of Heyting algebras with the top element as designated element. 

It is well known that a logical matrix $\langle \A, F\rangle$ can be regarded as a first-order structure, namely, as an algebra equipped with the interpretation of a predicate symbol $P(x)$. The intuitive reading of logical matrices suggests that the set of designated elements $F$ represents truth inside the set of truth-values $A$. Accordingly $P(x)$ can be understood as a \textit{truth predicate} and $F$ as the \textit{truth set} of $\langle \A, F\rangle$. Keeping this in mind, it makes sense to refer to the truth sets of a class of matrices.

One of the most striking achievements in the field of AAL is the discovery of the importance of the so-called \textit{Leibniz operator}. This is the map $\leibniz^{\A}\colon \mathcal{P}(A) \to \Con\A$, defined for every algebra $\A$, that sends every subset $F \subseteq A$ to the largest congruence $\theta$ of $\A$ such that $F$ is a union of blocks of $\theta$. The importance of the Leibniz operator comes from the fact that its behaviour on the deductive filters of a given logic $\LL$ determines interesting facts about the definability of \textit{logical equivalence} and of the \textit{truth sets} of the matrix semantics $\ModstarL$. This discovery led to the development of the so-called \textit{Leibniz hierarchy} \cite{Cz01,AAL-AIT-f,JM17book,JGRa11}, where logics are classified according to definability properties related to the behaviour of the Leibniz operator. Then the goal of this paper is to contribute to the systematic study of the aspects Leibniz hierarchy related to the definability of the truth sets of the matrix semantics $\ModstarL$. This program has already been considered in the AAL literature, especially in \cite{BlRe03,CzJa00,He93a,Ra06a}.

The starting point of our discussion is the following definability condition: we say that truth is \textit{almost parametrically equationally definable} in $\ModstarL$ if there is a set of equations $\btau(x, \vec{y})$ such that the \textit{non-empty} truth sets in $\ModstarL$ are exactly the sets of solutions of the equations $\btau(x, \vec{y})$ once we bound the parameters $\vec{y}$ by a universal quantifier. Similarly, we say that truth is \textit{equationally definable} in $\ModstarL$ when there is a set of equations $\btau(x)$ such that the truth sets in $\ModstarL$ are exactly the sets of solutions of the equations $\btau(x)$. It is clear that equational definability implies almost parametrized equational definability, and we show that the converse does not hold in general. The reader may wonder why we restrict the definition of almost parametrized equational definability to non-empty truth sets. This is because we prove that, when applied to all (possibly empty) truth sets of $\ModstarL$, the notion of equational and parametrized equational definability coincide. In particular, this implies that these two definability conditions are equivalent for logics with at least one tautology (Corollary \ref{Cor:Chap2:NoTheorem}).

Logics whose truth sets are equationally or almost parametrically equationally definable can be characterized by means of the behaviour of the Leibniz operator. More precisely, it turns out that truth is equationally (resp. almost parametrically equationally) definable in $\ModstarL$ if and only if the Leibniz operator is completely order-reflecting on (resp. on non-empty) deductive filters of $\LL$ over every algebra. This condition on the Leibniz operator can be equivalently restricted to the \textit{theories} of the logic $\LL$, i.e. to the filters of $\LL$ over the countably generated term algebra (Theorems \ref{Thm:Chap2:LeibnizUniversal} and \ref{Thm:Chap2:LeibnizEquational}). These results were first discovered by Raftery for equational definability in \cite{Ra06a}, where parametrized equational definability is not taken into account.

Until now we focused on logics whose truth sets can be defined by means of some linguistic 
translation of formulas into equations. This idea presents some analogy with the one of 
\textit{explicit definability} in first-order logic, in the sense that it requires that the definition 
of the truth sets is witnessed by some linguistic construction, i.e., by sets of equations. Now, Beth's definability theorem states that in first-order logic explicit definability and implicit definability coincide. Building on this analogy, it is natural to consider some suitable version of the notion of \textit{implicit definability} in the 
framework of truth sets of matrix semantics and to ask under which conditions these two kinds of definability 
coincide. Accordingly, given a logic $\LL$, we say that truth is \textit{implicitly definable} in $\ModstarL$ if the matrices in $\ModstarL$ are determined by their algebraic reduct. The analogy with Beth's definability theorem culminated in the discovery \cite{CzJa00,He93a,He97,Ra06a} that the notions of implicit and equational definability coincide when $\LL$ is a protoalgebraic logic (Theorem \ref{Thm:Chap3:Beth}).

Now, implicit definability can be characterized in terms of the behaviour of the Leibniz operator. More precisely, it has long been known that truth is implicitly definable in $\ModstarL$ if and only if the Leibniz operator is injective on the deductive filters of $\LL$ over \textit{every} algebra (Lemma \ref{Lem:Chap3:Implicit}). This fact posed the problem of whether the 
injectivity of the Leibniz operator transfers from the theories of a logic to its filters over 
arbitrary algebras \cite[Problem 1]{Ra06a}. The feeling that this question could have a positive answer was motivated by the fact that the main conditions on the Leibniz operator considered in the literature transfer from theories to filters over arbitrary algebras. In fact Czelakowski and Jansana provided 
in \cite{CzJa00} a positive answer to this problem, under the assumption of protoalgebraicity (Theorem \ref{Thm:Chap3:ProtoalgebraicLogics}). We solve this problem by showing that its answer depends on the cardinality of the language in which the logic is formulated. More precisely, if the language is countable, then the injectivity of the Leibniz operator transfers from theories to arbitrary filters (Theorem \ref{Thm:Chap3:TransferCountable}). On the other hand, we show that it is possible to construct counterexamples for logics expressed in uncountable languages (Section \ref{Sec:Failure}).

An intermediate definability condition that we take into account is the following: we say that truth is \textit{small} in $\ModstarL$ when the truth sets in this class are the smallest deductive filters of the logic $\LL$. We prove that this condition is equivalent to the fact that the Leibniz operator is  order-reflecting over deductive filters of every algebra (Lemma \ref{Lem:Chap3:Reflecting}). As it was the case for injectivity, the order-reflection of the Leibniz operator transfers from theories to filters over arbitrary algebras for logics expressed in a countable language (Theorem \ref{Thm:Chap3:TransferReflecting}), while there are counterexamples among logics whose language is uncountable. The work described until now originates an expansion of the Leibniz hierarchy with additional, weaker classes of logics corresponding to the definability conditions on the truth sets considered here. The expanded hierarchy is depicted in Figure \ref{Fig:LeibnizHierarchy2}.

\section{Matrix semantics}

Here we present a brief survey of the main concepts of abstract algebraic logic we will make use of along the article; a systematic exposition can be found for example in \cite{BP86,BP89,BP92,Cz01,AAL-AIT-f,FJa09,FJaP03b}. Fixed an algebraic type $\mathscr{L}$ and a \textit{countable} set $\Var$ of variables $x, y, z$, etc., we denote  $Fm$ the set of formulas over $\mathscr{L}$ built up with the variables $\Var$ and by $\Fm$ the corresponding absolutely free algebra. Moreover, given a formula $\varphi \in Fm$, we write $\varphi(x, \vec{z})$ if the variables of $\varphi$ are among $x$ and $\vec{z}$ and $x$ does not appear in $\vec{z}$. Keeping this in mind that, formally speaking, equations are pairs of formulas, we set $Eq \coloneqq Fm \times Fm$. From now on we assume that we are working with a fixed algebraic type.

By a \textit{logic} $\LL$ we understand a closure operator $C_{\lo}\colon \mathcal{P}(Fm) \to \mathcal{P}(Fm)$ which is \textit{structural} in the sense that $\sigma C_{\lo} \subseteq C_{\lo} \sigma$ for every endomorphism (or, equivalently, substitution) $\sigma \colon \Fm \to \Fm$. Given $\Gamma \cup \{\varphi\}\subseteq Fm$ sometimes we write $\Gamma \vdashL  \varphi$ instead of $\varphi \in C_{\lo}(\Gamma)$. Moreover, given $\Gamma \cup \{ \varphi, \psi \} \subseteq Fm$, we denote by $\Gamma, \varphi \sineq_{\LL} \psi, \Gamma$ the fact that both $\Gamma, \varphi \vdash_{\LL} \psi$ and $\Gamma, \psi \vdash_{\LL} \varphi$ are true. A formula $\varphi$ is a \textit{theorem} of $\LL$ if $\emptyset \vdash_{\LL} \varphi$. A logic $\LL$ is \textit{purely inferential} if it has no theorems. Given two logics $\lo$ and $\lo'$, we write $\lo \leq \lo'$ if $C_{\lo}(\Gamma) \subseteq C_{\lo'}(\Gamma)$ for every $\Gamma \subseteq Fm$. Since $\lo$ always denotes an arbitrary logic, we skip, in the formulation of our results, assumptions like ``let $\lo$ be a logic''.

We denote algebras with bold capital letters $\A$, $\B$, $\C$, etc.\ (with universes $A$, $B$, $C$, etc.\ 
respectively). The trivial algebra is denoted by $\boldsymbol{1}$. Given a logic $\lo$ and an algebra $\A$, we 
say that a set $F \subseteq A$ is a \textit{deductive filter} of $\lo$ over $\A$ when
\begin{gather*}
\text{if }\Gamma \vdashL  \varphi\text{, then for every homomorphism }h\colon \Fm\to \A,\\ 
\text{if }h[\Gamma] \subseteq F\text{, then }h(\varphi)\in F
\end{gather*}
for every $\Gamma \cup \{\varphi\} \subseteq Fm$. We denote by $\FFiLA$ the set of deductive filters of $\lo$ over $\A$, which turns out to be a closure system and, therefore, a complete lattice when ordered under set-theoretical inclusion. We denote by $\FiLA(\cdot)$ the closure operator of $\LL$-filter generation over the algebra $\A$. The filters over the algebra $\Fm$ are called \textit{theories} and their collection is denoted by $\Th\lo$ (instead of $\FFiL\Fm$). 

Given an algebra $\A$, we denote its congruence lattice by $\Con\A$. The identity relation on $\A$ is denoted by $\textup{Id}_{\A}$. A congruence $\theta \in \Con\A$ is \textit{compatible} with a set $F \subseteq A$ if for every $a, b\in A$
\[
\text{if }\langle a, b \rangle \in \theta \text{ and }a \in F\text{, then }b \in F.
\]
Given $F \subseteq A$, there exists always that largest congruence on $\A$ compatible with $F$. We 
 denote it by $\leibniz^{\A} F$ and call it the \textit{Leibniz congruence} of $F$ on 
$\A$. This notion naturally gives rise to a map $\leibniz^{\A} \colon \mathcal{P}(A) \to \Con \A$, called the \textit{Leibniz operator} which plays a fundamental role in this paper. 
\begin{Lemma}\label{Lem:Chap1:Inverse}
Let $f \colon \A \to \B$ be a homomorphism and $F \subseteq B$.  The following conditions hold:
\benormal
\item $f^{-1}\leibniz^{\B} F \subseteq \leibniz^{\A} f^{-1}[F]$. 
\item If $f$ is surjective, then $f^{-1}\leibniz^{\B} F = \leibniz^{\A} f^{-1}[F]$.
\enormal
\end{Lemma}

Given an algebra $\A$, a logic $\LL$, and a set $F \subseteq A$, the congruence
\[
\tarski^{\A}_{\LL}F \coloneqq \bigcap \{ \leibniz^{\A} G :  G \in \FFiLA \text{ and }F \subseteq G \}
\]
is called the \textit{Suszko congruence} of $F$ on $\A$ (relative to $\LL$). From the definition of the Suszko congruence, it follows that
\begin{equation}\label{Eq:SuszkoSmaller}
\tarski^{\A}_{\LL} F \subseteq \leibniz^{\A}F\text{ for every }F \in \FFiLA.
\end{equation}
The Suszko congruence is monotone in the sense that for every $F, G \subseteq A$ we have:
\begin{equation}\label{Eq:SuszkoMonotone}
\text{if }F \subseteq G\text{, then }\tarski^{\A}_{\LL} F \subseteq \tarski^{\A}_{\LL} G
\end{equation}
The Leibniz and Suszko congruences can be characterized in terms of the indiscernibility of elements with respect to filters in the following way. Given an algebra $\A$, a function $p \colon A^{n} \to A$ is a \index{polynomial function}\textit{polynomial function} if there are a natural number $m$, a formula $\varphi(x_{1}, \dots, x_{n+m})$ and elements $b_{1}, \dots, b_{m} \in A$ such that
\begin{align*}
p(a_{1}, \dots, a_{n}) = \varphi^{\A}(a_{1}, \dots, a_{n}, b_{1}, \dots, b_{m})
\end{align*}
for every $a_{1}, \dots, a_{n} \in A$. Observe that the notation $\varphi(x_{1},  \dots, x_{n+m})$ means just that the variables really occurring in $\varphi$ are among, but not necessarily all, $\{x_{1}, \dots, x_{n+m}\}$.
\begin{Lemma}\label{Lem:Polynomial} Let $\A$ be an algebra, $F \subseteq A$, and $a, b\in A$.
\benormal
\item $\langle a, b \rangle \in \leibniz^{\A}F\Longleftrightarrow (\?p(a)\in F\textrm{ if and only if }p(b)\in F\?)$ for every unary polynomial function $p$ of $\A$.
\item  $\langle a, b \rangle \in \tarski_{\LL}^{\A}F\Longleftrightarrow \FiLA(F \cup \{p(a)\}) = \FiLA(F \cup \{p(b)\})$ for every unary polynomial function $p$ of $\A$.
\enormal
\end{Lemma}

A pair $\langle \A, F \rangle$ is a \textit{matrix} if $F \subseteq A$. The \textit{reduction} of a matrix $\langle \A, F \rangle$ is the matrix $\langle \A, F \rangle^{\boldstar}\coloneqq \langle \A / \leibniz^{\A}F, F / \leibniz^{\A}F\rangle$. Moreover, a matrix $\langle \A, F \rangle$ is \textit{reduced} if $\leibniz^{\A}F = \textup{Id}_{\A}$. It can be easily 
seen that the reduction of a matrix is always reduced. The \textit{reduced models}, the \textit{Lindenbaum-Tarski models}, and the \textit{Suszko-reduced models} of a given a logic $\LL$ are, respectively, the following classes of matrices:
\begin{align*}
\ModstarL &\coloneqq \{ \langle \A, F \rangle : F \in \FFiLA \text{ and }\leibniz^{\A} F = \textup{Id}_{\A} \}\\
\LModstarL & \coloneqq \{ \langle \Fm / \leibniz \Gamma, \Gamma/ \leibniz \Gamma \rangle : \Gamma \in \Th\LL \}\\
\ModSL &\coloneqq \{ \langle \A, F \rangle : F \in \FFiLA \text{ and }\tarski^{\A}_{\LL} F = \textup{Id}_{\A} \}.
\end{align*}
The logic $\LL$ is complete with respect to any of the above classes of matrices.

\section{Definability with parameters}\label{Section Uni}

Before beginning, let us introduce some terminological convention, which will considerably simplify 
the formulation of the main results. We say that the Leibniz operator $\leibniz^{\A} \colon \mathcal{P}(A) \to \Con\A$ enjoys a certain set- or 
order-theoretic property over $\FFiLA$, if its restriction to $\FFiLA$ enjoys it. Moreover, we say that $\leibniz^{\A}$  \textit{almost} enjoys that 
property over $\FFiLA$ if its restriction to $\FFiLA \smallsetminus \{ \emptyset \}$ enjoys it. For example we will say that $\leibniz^{\A}$ is \textit{almost injective}\index{operator!(almost) injective} over $\FFiLA$ if $\leibniz^{\A}\colon  \FFiLA \smallsetminus \{ \emptyset \} \to \Con\A$ is injective. The reader may wonder why do we care so 
much about the empty filter. This is because we will be concerned with several examples of \textit{purely inferential} logics, i.e., logics without theorems. And 
it is easy to prove that $\emptyset \in \FFiLA$ if and only if $\LL$ is\index{logic!purely inferential} purely inferential. Therefore it will be often the case that the collection of 
deductive filters of our logics contains the empty-set, which represents a limit case and shall be eliminated in the formulation of the main results (that would be false otherwise). An 
analogous expedient will apply to matrices as follows. We say that a matrix $\langle \A, F\rangle$ is\index{matrix!almost trivial} \textit{almost trivial} if $F = \emptyset$. Accordingly, a matrix $\langle \A, F \rangle$ is \textit{non-almost trivial} when it is not almost trivial. Observe that $\langle \boldsymbol{1}, \emptyset \rangle$ is the unique matrix that is both reduced and almost trivial. A class of matrices $\class{M}$ \textit{almost} enjoys a certain 
property, if every non-almost trivial member of $\class{M}$ enjoys it.

\begin{Definition} A \textit{parametrized equational translation} is a set $\btau (x, \vec{y}) \subseteq Eq$ of equations in a distinguished variable $x$ with parameters $\vec{y}$. An \textit{equational translation} is a parametrized equational translation without parameters $\vec{y}$.
\end{Definition}
Parametrized equational translations witness the definability of truth sets in classes of matrices by bounding parameters (if any) by an universal quantifier and considering the solutions of the resulting universally quantified equations. More precisely, given a parametrized equational translation $\btau(x, \vec{y})$ and an algebra $\A$, we let
\begin{equation}\label{Solutions}
\btau(\A) \coloneqq \{ a\in A : \A \vDash \btau(a, \vec{c}) \textrm{ for every }\vec{c}\in A\}.
\end{equation}
Observe that if $\btau(x, \vec{y}) = \btau(x)$ is an equational translation, then (\ref{Solutions}) simplifies to the following:
\begin{align*}
\btau(\A) = \{ a\in A : \A \vDash \btau(a) \}.
\end{align*}

\begin{Definition}\label{Def:Chap2:UnivDedf} A parametrized equational translation (resp.\ equational translation) $\btau$ \textit{defines truth} in a class of matrices $\class{M}$, if $\btau(\A) = F$ for every $\langle \A, F\rangle \in \class{M}$. Truth is \textit{parametrically equationally \textup{(resp.} equationally\textup{)} definable} in $\class{M}$ if there is a parametrized equational translation (resp.\ equational translation) that  defines truth in $\class{M}$.
\end{Definition}

\begin{exa}[\textsf{Lattices}]\label{Ex:LatticeFirstExample}
Let $\A$ be a lattice with a maximum element $a$. Then consider the matrix $\langle \A, \{ a \} \rangle$. For every $b \in A$ we have that
\begin{eqnarray*}
b = a &\Longleftrightarrow & c \leq b \text{ for every }c \in A\\
& \Longleftrightarrow & c \land b = c  \text{ for every }c \in A\\
& \Longleftrightarrow & \A \vDash \btau(b, c) \text{ for every }c \in A
\end{eqnarray*}
where $\btau(x, \vec{y})$ is the parametrized equational translation $\{ x \land y \thickapprox y \}$. This shows that $\btau$ defines truth in $\langle \A, \{ a \} \rangle$. On the other hand if $\A$ is non-trivial, there is no equational translation that defines truth in $\langle \A, \{ a \}\rangle$. This is due to the fact that (up to equivalence) the unique lattice equation in variable $x$ is $x \thickapprox x$. The situation changes if we add a constant $1$ to the type of $\A$. In particular, let $\A^{+}$ be the expansion of $\A$ where $1$ is interpreted as $a$. Then truth is equationally definable in $\langle \A^{+}, \{ a \} \rangle$ by the equational translation $\btau(x) = \{ x \thickapprox 1 \}$.
\qed
\end{exa}

Observe that when truth is almost parametrically equationally definable in $\class{M}$, the non-almost trivial matrices in $\class{M}$ are determined by their algebraic reduct. More precisely, if $\langle \A, F\rangle, \langle \A, G\rangle \in \class{M}$ are non-almost trivial, then $F = G$. This observation will be used in several proofs. It is clear that if truth is equationally definable in $\class{M}$, then it is parametrically equationally definable too. 
 
We will be 
interested in logics $\LL$ for which truth is parametrically equationally or equationally definable in $\ModstarL$. For this reason it will be convenient to introduce some terminological  convention. We say that the \textit{truth sets} of a logic $\LL$ are \textit{parametrically equationally} (resp. equationally) \textit{definable}, as an abbreviation for the fact that truth is parametrically equationally (resp. equationally) definable in $\ModstarL$. In the case of equational definability, logics that satisfy this property have been studied in depth in the literature \cite{CzJa00,He93a,He93,Ra06a}.\footnote{In \cite{Ra06a} Raftery calls \textit{truth-equational} the logics $\LL$ for which truth is equationally definable in $\ModstarL$. Here we prefer not to give them any particular name, in order to obtain a more uniform naming scheme when dealing with different kinds of definability conditions.} One of the main outcomes of this research line was the discovery of Raftery \cite[Theorem 28]{Ra06a} that logics whose truth sets are equationally definable can be characterized in terms of the behaviour of the Leibniz operator as follows:
\begin{Definition}
Let $\boldsymbol{X}$ and $\boldsymbol{Y}$ be complete lattices and $f \colon X \to Y$ be a map.  $f$ is \textit{completely order-reflecting} if for every $A \cup \{ b \} \subseteq X$,
\[
\text{if }\bigwedge_{a \in A} f(a)\leq f(b)\text{, then }\bigwedge_{a \in A } a \leq b.
\]
\end{Definition}

\begin{Theorem}[Raftery]\label{Thm:Chap2:LeibnizEquational} The following conditions are equivalent:

\benroman
\item The truth sets of $\LL$ are equationally definable.
\item $\leibniz^{\A} \colon \FFiLA \to \Con\A$ is completely order-reflecting, for every algebra $\A$.
\item $\leibniz \colon \Th\LL \to \Con\Fm$ is completely order-reflecting.
\end{enumerate}
In this case $\btau(x) \coloneqq \sigma_{x}\tarski_{\LL} C_{\LL}\{x\}$ defines truth in $\ModstarL$, where $\sigma_{x}$ is the substitutions sending every variable to $x$.
\end{Theorem}

The main goal of this section will be to provide an analogous characterization of logics whose truth sets are almost parametrically equationally definable (Theorem \ref{Thm:Chap2:LeibnizUniversal}). One may wonder why we are interested in logics $\LL$, whose truth sets are \textit{almost} parametrically equationally definable, and not simply parametrically equationally definable. This is because the notions of equational and parametrized equational definability coincide, when they are applied to the truth sets of a logic (Corollary \ref{Cor:Chap2:NoTheorem}).

Formally speaking an equation equation $\epsilon \thickapprox \delta$ is just a pair $\langle \epsilon, \delta \rangle$. Thus, given an algebra $\A$ and a tuple $\vec{a} \in A$, the expression $\epsilon^{\A}(\vec{a}) \thickapprox \delta^{\A}(\vec{a})$ will denote the pair $\langle \epsilon^{\A}(\vec{a}), \delta^{\A}(\vec{a}) \rangle \in A \times A$. Keeping this in mind, we have the following:
\begin{Lemma}\label{Lem:Chap2:technical} Let $\btau(x, \vec{y})$ be a parametrized equational translation and $\langle \A, F\rangle$ a non-almost trivial matrix. $\btau(x, \vec{y})$ defines truth in the reduction $\langle \A, F\rangle^{\boldstar}$ if and only if for every $a \in A$,
\[
a\in F \Longleftrightarrow \btau^{\A}(a, \vec{c}) \subseteq \leibniz^{\A} F\text{ for every }\vec{c}\in A.
\]
\end{Lemma}
\begin{proof}
Apply the fact that $\leibniz^{\A}F$ is compatible with $F$.
\end{proof}

The first fact that it is worth to remark is the following:

\begin{Lemma}\label{Lem:Chap2:Transfer} A parametrized equational translation almost defines truth in $\ModstarL$ if and only if it almost defines truth in $\LModstarL$.
\end{Lemma}
\begin{proof}
It will be enough to check the  ``if'' part.  Let $\btau$ be a parametrized equational 
translation which almost defines truth in $\LModstarL$. We want to show that 
$\btau$ almost defines truth in $\ModstarL$ too. By Lemma \ref{Lem:Chap2:technical} this amounts to proving that for every algebra $\A$, $F \in \FFiL\A \smallsetminus \{ \emptyset \}$ and $a\in A$:
\begin{equation}\label{Eq:Chap2:za}
a\in F \Longleftrightarrow\btau^{\A}(a, \vec{c})\subseteq \leibniz^{\A}F\textrm{ for every }\vec{c}\in A.
\end{equation}
First recall that $\LModstarL$ is (up to isomorphism) the class of countably generated reduced models of $\LL$. Together with Lemma \ref{Lem:Chap2:technical}, this implies that (\ref{Eq:Chap2:za}) holds in case $\A$ is countably generated.

Then consider the case where $\A$ is not countably generated. We begin by proving the ``if'' part of (\ref{Eq:Chap2:za}). Let $a\in A$ and $F \in \FFiL\A \smallsetminus \{ \emptyset \}$ be such that $\btau^{\A}(a, \vec{c})\subseteq \leibniz^{\A} F$ for every $\vec{c}\in A$. Then choose any $b\in F$ and consider the subalgebra of $\B$ of $\A$ generated by $\{a, b\}$. Let $G\coloneqq F\cap B$. Observe that $G \in \FFiL \B \smallsetminus \{ \emptyset \}$, since $b \in G$. Moreover, we have $\btau^{\B}(a, \vec{c})\subseteq \leibniz^{\A} F \cap (B\times B) \subseteq \leibniz^{\B} G$ for every $\vec{c} \in B$. Since $\B$ is countably generated, we conclude that $a \in G \subseteq F$.

Now we prove the ``only if'' part of (\ref{Eq:Chap2:za}). Suppose that $a\in F$. Then consider 
$\epsilon \thickapprox \delta \in \btau$ and $\vec{c}\in A$. Let also $p(x)\colon \A\to \A$ be a unary polynomial 
function. By definition there is an $(n+1)$-ary term $\varphi$ and a   
sequence $\vec{e}$ of $n$ elements of $A$ such that $\varphi^{\A}(b, \vec{e})= p(b)$ for every $b \in A$. We will prove that
\begin{equation}\label{Eq:Chap2:zb}
p(\epsilon^{\A}(a, \vec{c}))\in F\Longleftrightarrow p(\delta^{\A}(a, \vec{c}))\in F.
\end{equation}
\noindent First suppose that $p(\epsilon^{\A}(a, \vec{c}))\in F$. Then consider the subalgebra $\B$ of $\A$ generated by 
$\{a, e_{1}, \dots, e_{n}, c_{1}, \dots, c_{k}\}$, where $c_{1}, \dots c_{k}$ are the 
elements of $\vec{c}$ corresponding to the variables in $\vec{y}$ occurring in $\epsilon$ and $\delta$. Then let $G\coloneqq F\cap B$. We know that $G\in \FFiL \B \smallsetminus \{ \emptyset \}$, because $a\in G$.  Since $\B$ is countably generated and $a \in G$, we have that
\[
\langle \epsilon^{\B}(a, c_{1}, \dots, c_{k}), \delta^{\B}(a, c_{1}, \dots, c_{k}) \rangle \in \leibniz^{\B}G.
\]
Since $p(\epsilon^{\B}(a, \vec{c}))\in G$, by compatibility we obtain that $p(\delta^{\B}(a, \vec{c})) \in G \subseteq F$. This establishes condition (\ref{Eq:Chap2:zb}). By point 1 of Lemma \ref{Lem:Polynomial}, we conclude that $\langle \epsilon^{\A}(a, \vec{c}), \delta^{\A}(a, \vec{c})\rangle \in \leibniz^{\A}F$.
\end{proof}

The following technical result is stated without a detailed proof in \cite[Proposition 1.5(8)]{Cze03} and will be needed in the sequel.

\begin{Lemma}[Czelakowski]\label{Lem:Chap1:ThecnicalSuszko} $\sigma \tarski_{\LL} C_{\LL}\{x\} \subseteq \tarski_{\LL} C_{\LL}\{\sigma x\}$ for every substitution $\sigma$.
\end{Lemma}
\begin{proof}
Consider a pair $\langle \varphi, \psi \rangle \in \tarski_{\LL} C_{\LL}\{x\}$. We have to prove that $\langle \sigma \varphi, \sigma \psi \rangle \in \tarski_{\LL} C_{\LL}\{\sigma x\}$. By point 2 of Lemma \ref{Lem:Polynomial} it will be enough to check that 
\[
\gamma(\sigma(\varphi), \vec{z}), \sigma x \sineq_{\LL} \sigma x, \gamma(\sigma(\psi), \vec{z}) \text{, for every }\gamma(x, \vec{z}) \in Fm.
\]
To this end, consider $\gamma(x, \vec{z}) \in Fm$ and a new substitution $\sigma$ such that:

\benormal
\item $\sigma$ and $\sigma'$ coincide on the variables actually occurring in $\varphi$ and $\psi$.
\item For every variable $v \ne x$ actually occurring in $\gamma$, there is a variable $u$ such that $\sigma'(u) = v$.
\enormal
Now, consider the formula $\delta$ obtained by replacing in $\gamma$ each variable $v \ne x$ by the corresponding $u$. Applying point 2 of Lemma \ref{Lem:Polynomial} to the fact that $\langle \varphi, \psi \rangle \in \tarski_{\LL} C_{\LL}\{x\}$, we obtain that
\[
x, \delta(\varphi, \vec{u}) \sineq_{\LL}\delta(\psi, \vec{u}), x.
\]
By structurality we obtain that
\[
\sigma'x, \sigma'\delta(\varphi, \vec{u}) \sineq_{\LL}\sigma'\delta(\psi, \vec{u}), \sigma'x.
\]
But this is exactly $\gamma(\sigma(\varphi), \vec{z}), \sigma x \sineq_{\LL} \sigma x, \gamma(\sigma(\psi), \vec{z})$.
\end{proof}

The next result provides a characterization of the logics whose truth sets are almost parametrically equationally definable in terms of the behaviour of the Leibniz operator. One can read the result as stating that the class of logics whose truth sets are almost parametrically equationally definable belongs to the Leibniz hierarchy.

\begin{Theorem}\label{Thm:Chap2:LeibnizUniversal} The following conditions are equivalent:

\benroman
\item The truth sets of $\LL$ are almost parametrically equationally definable.
\item $\leibniz^{\A} \colon \FFiLA \to \Con\A$ is almost completely order-reflecting, for every algebra $\A$.
\item $\leibniz \colon \Th\LL \to \Con\Fm$ is almost completely order-reflecting.
\end{enumerate}
In this case $\btau(x, \vec{y}) \coloneqq \tarski_{\LL} C_{\LL}\{x\}$ almost defines truth in $\ModstarL$.
\end{Theorem}
\begin{proof}
(i)$\Rightarrow$(ii): Consider an arbitrary algebra $\A$ and let $\mathcal{F} \cup \{G\} \subseteq \FFiLA\smallsetminus \{\emptyset\}$ such that $\bigcap \{ \leibniz^{\A} F : F\in \mathcal{F} \} \subseteq \leibniz^{\A} G$. We proceed to show that $\bigcap \mathcal{F} \subseteq G$. To this end, consider $a\in \bigcap \mathcal{F}$. From Lemma \ref{Lem:Chap2:technical} and the assumptions it follows that
\[
\btau^{\A}(a, \vec{c})\subseteq \leibniz^{\A} F\text{ for every }\vec{c}\in A\text{ and }F\in \mathcal{F}.
\]
This implies that $\btau^{\A}(a, \vec{c})\subseteq \leibniz^{\A} G$ for every $\vec{c}\in A$. With another application of Lemma \ref{Lem:Chap2:technical} we conclude that $a\in G$.

(ii)$\Rightarrow$(iii): Straightforward. (iii)$\Rightarrow$(i): Choose a variable $x$ and define $\btau(x, \vec{y}) \coloneqq \tarski_{\LL} C_{\mathcal{L}}\{x\}$ where $\vec{y}$ is the list of all variables different from $x$. Thanks to Lemma \ref{Lem:Chap2:Transfer} it will be enough to prove that $\btau$ almost defines truth in $\LModstarL$. By Lemma \ref{Lem:Chap2:technical} this reduces to proving the following:
\begin{equation}\label{Eq:Chap2:Leibniz}
\Gamma \vdash_{\LL} \varphi \textrm{ if and only if }\btau(\varphi, \vec{\gamma}) \subseteq \leibniz\Gamma \textrm{ for every } \vec{\gamma} \in Fm
\end{equation}
for every $\varphi\in Fm$ and $\Gamma \in \Th\LL \smallsetminus \{\emptyset\}$.
For the ``only if'' part of (\ref{Eq:Chap2:Leibniz}) suppose that $\Gamma \vdash_{\LL} \varphi$. Then consider any sequence $\vec{\gamma} \in Fm$ and let $\sigma$ be a substitution sending $x$ to $\varphi$ and $\vec{y}$ to $\vec{\gamma}$. Applying Lemma \ref{Lem:Chap1:ThecnicalSuszko}, we obtain that 
\begin{align*}
\btau(\varphi, \vec{\gamma}) &= \sigma \btau(x, \vec{y}) = \sigma \tarski_{\LL} C_{\LL}\{ x \}\subseteq \tarski_{\LL} C_{\LL}\{\sigma x\}\\
&  = \tarski_{\LL} C_{\LL}\{\varphi\} \subseteq \tarski_{\LL} \Gamma \subseteq \leibniz \Gamma.
\end{align*}

Then we turn to prove the ``if'' part of (\ref{Eq:Chap2:Leibniz}). Suppose that $\btau(\varphi, \vec{\gamma}) \subseteq \leibniz\Gamma$ for every sequence $\vec{\gamma} \in Fm$. Recall that $\Gamma \ne \emptyset$. Then we can choose a formula $\psi\in \Gamma$ and consider the substitution $\sigma$ defined as
\begin{align*}
\sigma(z) = \left\{ \begin{array}{ll}
\varphi & \textrm{if $z=x$}\\
\psi & \textrm{otherwise}\\
\end{array} \right.
\end{align*}
for every variable $z$. From Lemma \ref{Lem:Chap1:Inverse} and the assumption it follows that
\[
\tarski C_{\mathcal{L}}\{x\} = \btau(x, \vec{y}) \subseteq \sigma^{-1} \leibniz\Gamma \subseteq \leibniz \sigma^{-1}\Gamma.
\]
Since inverse images of theories under substitutions are theories, we know that $\sigma^{-1}\Gamma \in \Th\LL$. Moreover, observe that $y\in \sigma^{-1}\Gamma$ for every variable different from $x$. Thus $\sigma^{-1}\Gamma\ne \emptyset$. Therefore we can apply the fact that $\leibniz$ is completely order-reflecting over $\Th\LL\smallsetminus \{\emptyset\}$ and get $C_{\mathcal{L}}\{x\} \subseteq \sigma^{-1}\Gamma$. This yields that $\varphi\in \Gamma$ and concludes the proof of condition (\ref{Eq:Chap2:Leibniz}).
\end{proof}

Combining Theorems \ref{Thm:Chap2:LeibnizUniversal} and \ref{Thm:Chap2:LeibnizEquational} we can prove a surprising result, namely that truth is parametrically equationally definable in the whole class of Leibniz-reduced models of a logic if and only if it is equationally definable in it.

\begin{Corollary}\label{Cor:Chap2:NoTheorem} The following conditions are equivalent:

\benroman
\item The truth sets of $\LL$ are equationally definable.
\item The truth sets of $\LL$ are parametrically equationally definable.
\item The truth sets of $\LL$ are almost parametrically equationally definable and $\LL$ has theorems.
\end{enumerate}
In particular, if truth is equationally definable in $\ModstarL$, then $\LL$ has theorems.
\end{Corollary}
\begin{proof}
(i)$\Rightarrow$(ii): Straightforward. (ii)$\Rightarrow$(iii): Suppose towards a contradiction that $\LL$ is purely inferential, i.e., that $\emptyset \in \Th\LL$. Then both $\langle \boldsymbol{1}, \{1\}\rangle$ and $\langle \boldsymbol{1}, \emptyset\rangle$ are reduced models of $\LL$. But this contradicts the fact that truth is parametrically equationally definable in $\ModstarL$. (iii)$\Rightarrow$(i): Together with Theorem \ref{Thm:Chap2:LeibnizUniversal}, the assumption implies that $\leibniz$ is completely order-reflecting over $\Th\LL$. Thus with an application of Theorem \ref{Thm:Chap2:LeibnizEquational} we are done.
\end{proof}

From Corollary \ref{Cor:Chap2:NoTheorem} it follows that the notion of parametrized equational definability makes sense only for \textit{purely inferential} logics since, in the presence of theorems, it collapses into that of equational definability. It is worth to observe that the class of all logics, whose truth sets are almost parametrically equationally definable, is the first known class in the Leibniz hierarchy that admits non-trivial purely inferential logics. More in detail, until to now the weakest conditions considered in the study of the Leibniz hierarchy were \textit{protoalgebraicity} and \textit{having truth sets equationally definable}. It is well known that the unique purely inferential protoalgebraic logic (in a given language) is the almost inconsistent one,\footnote{Within a fixed algebraic language, the \textit{almost inconsistent logic} is the unique logic $\LL$ such that $\Th\LL = \{ \emptyset, Fm \}$.} while the fact that logics whose truth sets are equationally definable have theorems is stated in Corollary \ref{Cor:Chap2:NoTheorem}.
 
\begin{Corollary}\label{Cor:Chap2:Conservative}
Let $\LL$ be a logic whose truth sets are almost parametrically equationally definable. There is a conservative expansion $\LL'$ of $\LL$, where the expansion consists in adding a new constant symbol $1$, such that truth is equationally definable in $\ModstarL'$.  
\end{Corollary}
\begin{proof}
Let $\btau(x, \vec{y})$ be the parametrized equational translation that almost defines truth in $\ModstarL$. Then let $\class{K}$ be the class of algebras obtained by expanding the algebras in $\AlgstarL$ with a fresh constant $1$, that is interpreted arbitrarily in the set of solutions of $\btau$. Observe that an algebra $\A \in \AlgstarL$ can be expanded in different ways if $\btau(\A)$ has more than one element. Then let $\LL'$ be the logic determined by the following class of matrices:
\[
\{ \langle \A, \btau(\A) \rangle : \A \in \class{K} \}.
\]
It is easy to see that $\LL'$ is a conservative expansion of $\LL$. Since the expansion consists in adding a constant symbol and the presence of constants does not affect congruences, we obtain that if $\langle \A, F\rangle \in \ModstarL'$, then $\langle \A', F\rangle \in \ModstarL$ where $\A'$ is the $1$-free reduct of $\A$. Thus from the fact that the truth sets of $\LL$ are almost parametrically equationally definable, we can assume that the same holds for $\LL'$. Now observe that $\LL'$ has theorems, since $\emptyset \vdash_{\LL'}1$. Thus, with an application of Corollary \ref{Cor:Chap2:NoTheorem} we conclude that truth is equationally definable in $\ModstarL'$.
\end{proof}

\begin{problem}
In \cite[Theorem 11]{Ra06a} it is shown that the truth sets of a logic $\LL$ are equationally definable if and only if the $\tarski_{\LL}^{\A} \colon \FFiLA \to \Con\A$ is injective for every algebra $\A$. Is it possible to produce a similar characterization of almost parametrized equational definability?
\end{problem}

\section{Semilattice-based examples}\label{Sec:UniExamples}

In this section we review a family of natural examples of logics whose truth sets are almost parametrically equationally, but not equationally, definable. In the light of Corollary \ref{Cor:Chap2:NoTheorem} we know that all these examples need to be purely inferential.

\begin{exa}[\textsf{Distributive Lattices}]\label{Ex:Chap2:DistributiveLattices} Let be the $\langle \land, \lor \rangle$-fragment of classical propositional logic. For every non-almost trivial matrix $\langle \A, F \rangle \in \Modstar \mathcal{CPC}_{\land \lor}$ the following conditions hold:

\benroman
\item $\A$ is a distributive lattice with a maximum $1$.
\item $F = \{1\}$.
\item For every $a, b\in A$, if $a<b$, then  there is $c\in A$ such that $a\lor c\ne 1$ and $a\lor c=1$.
\end{enumerate}
This was proved in \cite[Pag.\ 127]{FGV91}, but see \cite{FV91} for further information on the logic $\mathcal{CPC}_{\land \lor}$. In particular, this result implies that the truth sets of $\mathcal{CPC}_{\land \lor}$ are almost parametrically equationally definable through the parametrized equational translation $\btau(x, \vec{y})=\{x \land y \thickapprox y\}$, as shown in Example \ref{Ex:LatticeFirstExample}.
\qed
\end{exa}

\begin{exa}[\textsf{Semilattices}]\label{Ex2} 
Let $\mathcal{CPC}_{\land}$ be the $\langle \land\rangle$-fragment of classical propositional logic. Let also $\boldsymbol{2} = \langle \{0, 1\}, \land \rangle$ be the two-element meet semilattice with $0 < 1$. Every non-almost trivial member of $\Modstar\mathcal{CPC}_{\land}$ is an isomorphic copy either of $\langle \mathbf{2}, \{ 1 \}\rangle$ or of $\langle \boldsymbol{1}, \{1\}\rangle$. This was first claimed in \cite[Pag. 68-69]{Rt93}, but see \cite[Corollary 6.3]{FMo14c} or \cite[Lemma 3.1]{Mo13a} for an explicit proof. It follows that the truth sets of $\mathcal{CPC}_{\land}$ are almost parametrically equationally definable by the parametrized equational translation $\btau(x, \vec{y})=\{x \land y \thickapprox y\}$.

Let\index{3@$\mathcal{CPC}_{\lor}$} $\mathcal{CPC}_{\lor}$ be the $\langle \lor\rangle$-fragment of classical propositional logic. The non-almost trivial members of $\Modstar\mathcal{CPC}_{\lor}$ can be characterized exactly as those of $\mathcal{CPC}_{\land\lor}$ in Example \ref{Ex:Chap2:DistributiveLattices}, but replacing \textit{distributive lattice} by \textit{semilattice} \cite[Pag. 68-69]{Rt93}. It follows that the truth sets of $\mathcal{CPC}_{\lor}$ are almost parametrically equationally definable by the parametrized equational translation $\btau(x, \vec{y})=\{x \lor y \thickapprox x \}$.
\qed
\end{exa}

Until now we saw that the truth sets of $\mathcal{CPC}_{\land \lor}, \mathcal{CPC}_{\land}$ and $\mathcal{CPC}_{\lor}$ are almost parametrically equationally  definable. The fact that the truth predicates of these logics are not equationally definable follows from the fact that (up to equivalence) the unique lattice or semilattice equation in variable $x$ is $x \thickapprox x$. The examples considered so far are particular instances of the following general phenomenon:

\begin{Lemma}\label{Lem:Chap2:JoinReduct}
Let $\class{K}$ be a class of algebras with meet (join) semilattice reduct. The truth sets of the logic $\LL$ determined by the class of matrices
\[
\{ \langle \A, \{ 1 \} \rangle : \A \in \class{K}\text{ has a top element }1 \}
\]
are almost parametrically equationally definable.
\end{Lemma}
\begin{proof}
We claim that $\AlgstarL \subseteq \VVV(\class{K})$. To prove this, consider an equation $\varphi \thickapprox \psi$ that holds in $\class{K}$. We want to show that for every $\vec{a} \in A$, we have $\langle \varphi^{\A}(\vec{a}), \psi^{\A}(\vec{a}) \rangle \in \textup{Id}_{\A} = \leibniz^{\A}F$. By Lemma \ref{Lem:Polynomial} this amounts to showing that for every unary polynomial function $p$ of $\A$ we have that
\begin{equation}\label{Eq:LeibnizApplication}
p(\varphi^{\A}(\vec{a})) \in F \Longleftrightarrow p(\psi^{\A}(\vec{a})) \in F.
\end{equation}
We know that there are a formula $\delta(x, \vec{z})$ and a tuple $\vec{c} \in A$ such that $p(x) = \delta^{\A}(x, \vec{c})$. We can assume w.l.o.g.\ that the tuple $\vec{z}$ does not contain any variable occurring in $\varphi$ or $\psi$. Since the equation $\varphi \thickapprox \psi$ holds in $\class{K}$, we have that
\begin{equation}\label{Eq:ValidityInK}
\class{K} \vDash \delta(\varphi, \vec{z}) \thickapprox \delta(\psi, \vec{z}).
\end{equation}
The fact that $\LL$ is defined by a class of matrices that have algebraic reducts in $\class{K}$, together with (\ref{Eq:ValidityInK}) implies that
\[
\delta(\varphi, \vec{z}) \sineq_{\LL} \delta(\psi, \vec{z}).
\]
This establishes condition (\ref{Eq:LeibnizApplication}) and concludes the proof of the claim.

Now observe that $\LL$ is determined by matrices whose filters are singletons. A standard argument \cite[Theorem 10]{AFRM15} shows that $F$ is a singleton for every non-almost trivial $\langle \A, F\rangle \in \ModstarL$. Then consider one of these models $\langle \A, F\rangle$. From our claim it follows that $\A \in \VVV(\class{K})$, thus $\A$ has a meet-semilattice reduct $\langle A, \land \rangle$. Observe that by definition of $\LL$
\[
x, y \vdash_{\LL} x \land y \qquad x \land y \vdash_{\LL} x \qquad x \land y \vdash_{\LL} y.
\]
Thus $F$ is a filter of $\langle A, \land \rangle$ which, moreover, is a singleton. Hence $\langle A, \land \rangle$ has a top element $1$ such that $F = \{ 1 \}$. We conclude that truth is almost parametrically equationally definable in $\ModstarL$ via $\btau(x, \vec{y}) \coloneqq \{ x \land y \thickapprox y \}$.
\end{proof}

It makes sense to wonder whether there are examples of meaningful logics whose truth set are almost parametrically equationally (but not equationally), except from those that fall under the scope of Lemma \ref{Lem:Chap2:JoinReduct}. One of them comes from the study of bilattices, i.e., algebras which have two lattice-theoretic order relations.

\begin{exa}[\textsf{Distributive Bilattices}]\label{bilattices} An\index{bilattice} algebra $\boldsymbol{L} = \langle L, \land, \lor, \otimes, \oplus\rangle$ is a \textit{pre-bilattice} if $\langle L, \land, \lor\rangle$ and $\langle L, \otimes, 	\oplus\rangle$ are lattices. In this case we will denote by $\leq$ the lattice order associated with $\langle L, \land, \lor\rangle$ and by $\sqsubseteq$ the one associated with $\langle L, \otimes, \oplus\rangle$. An algebra $\A = \langle A, \land, \lor, \otimes, \oplus, \lnot\rangle$ is a \textit{bilattice} if $\langle A, \land, \lor, \otimes, \oplus\rangle$ is a pre-bilattice such that $\lnot \lnot a = a$ and
\[
a \leq b \Longrightarrow (\lnot b\leq \lnot a \text{ and } \lnot a \sqsubseteq\lnot b)
\]
for every $a, b \in A$ \cite{Gi88}. A bilattice is \textit{distributive} if the four lattice operations satisfy all the combined distributive axioms. Distributive bilattices form a variety which we denote by $\class{DBL}$. This variety is generated by the bilattice $\B_{4}$ with universe $\{0, 1, a, b\}$, where $\langle B_{4}, \sqsubseteq\rangle$ is the four-element diamond bounded by $0<1$, while $\langle B_{4}, \leq \rangle$ is the four-element diamond bounded by $a<b$, and $\lnot$ interchanges $a$ and $b$ and is the identity on $\{ 0, 1 \}$.

The so-called \textit{logic of distributive bilattices}\index{9@$\mathcal{LB}$} $\mathcal{LB}$ is defined through the matrix $\langle \B_{4}, \{ b, 1 \} \rangle$, see \cite{BoRi11,Riv-PhD}. By Corollary \ref{Cor:Chap2:NoTheorem} truth is not equationally definable in $\Modstar\mathcal{LB}$, since $\mathcal{LB}$ is purely inferential. This is a consequence of the fact that $\{ 0 \}$ is the universe of a subalgebra of $\B_{4}$. Our goal will be to prove that truth is almost parametrically equationally definable in $\Modstar\mathcal{LB}$. To this end, we need to recall the following construction from \cite{BoRi11}. Given a distributive lattice $\boldsymbol{L} = \langle L, \sqcap, \sqcup \rangle$, let $\boldsymbol{L} \odot \boldsymbol{L} = \langle L \times L, \land, \lor, \otimes, \oplus, \lnot\rangle$ be the twist structure defined as
\begin{align*}
\langle a_{1}, a_{2}\rangle \land \langle b_{1}, b_{2}\rangle & \coloneqq \langle a_{1} \sqcap b_{1}, a_{2}  \sqcup b_{2}\rangle\\
\langle a_{1}, a_{2}\rangle \lor \langle b_{1}, b_{2}\rangle &\coloneqq \langle a_{1}  \sqcup b_{1}, a_{2} \sqcap b_{2}\rangle\\
\langle a_{1}, a_{2}\rangle \otimes \langle b_{1}, b_{2}\rangle &\coloneqq \langle a_{1} \sqcap b_{1}, a_{2} \sqcap b_{2}\rangle\\
\langle a_{1}, a_{2}\rangle \oplus \langle b_{1}, b_{2}\rangle &\coloneqq \langle a_{1} \sqcup b_{1}, a_{2} \sqcup b_{2}\rangle\\
\lnot \langle a_{1}, a_{2}\rangle &\coloneqq \langle a_{2}, a_{1}\rangle
\end{align*}
for every $\langle a_{1}, a_{2}\rangle, \langle b_{1}, b_{2}\rangle \in L \times L$. It turns out that $\boldsymbol{L} \odot \boldsymbol{L} \in \class{DBL}$.

It turns out that the parametrized equational translation $\btau(x, \vec{y}) \coloneqq \{ (x \oplus y) \land x \thickapprox x \oplus y \}$ almost defines  truth in $\Modstar\mathcal{LB}$. To see this, observe that if $\langle \A, F \rangle \in \Modstar\mathcal{LB}$ is non-almost trivial, then $\A \cong \boldsymbol{L} \odot \boldsymbol{L}$ for some lattice $\boldsymbol{L}$ such that the following conditions hold:

\benroman
\item $\boldsymbol{L}$ is a distributive lattice with maximum $1$.
\item $F\cong \{1\} \times L$.
\item For every $a, b\in L$, if $a < b$, then there is $c\in L$ such that $a \sqcup c < b \sqcup c = 1$.
\end{enumerate}
This was proved in \cite[Theorem 4.13]{BoRi11}. For the sake of simplicity, we assume that $\A = \boldsymbol{L}\odot \boldsymbol{L}$. For every $\langle a_{1}, a_{2}, \rangle , \langle b_{1}, b_{2}\rangle \in A$ we have that
\[
\A \vDash \btau(\langle a_{1}, a_{2} \rangle,\langle b_{1}, b_{2}\rangle) \Longleftrightarrow b_{1} \leq a_{1}.
\]
This yields that $\btau(\A) = \{1\} \times L$, matching condition (ii) above.
\qed
\end{exa}

All the logics considered so far were equipped either with a disjunction or with a conjunction, which was interpreted as a semilattice operation. Remarkably, within the landscape of Fregean logics, the presence of a weak disjunction (conjunction) forces the truth sets to be almost parametrically equationally definable. Recall that a logic $\LL$ is \textit{Fregean} if for every set of formulas $\Gamma$ the following relation is a congruence of the term algebra $\Fm$:
\[
\frege_{\LL} \Gamma \coloneqq \{ \langle \varphi, \psi\rangle \in Fm^{2} : \Gamma, \varphi \sineq_{\LL} \psi, \Gamma \}.
\]
Equivalently, a logic is Fregean when $\tarski_{\LL}\Gamma = \frege_{\LL} \Gamma$ for every $\Gamma \in \Th\LL$. Fregean logics have been studied in depth \cite{Ba03,Cz01,CzPi04,CzePi04}, but here we will rely only on their definition.

\begin{Definition}
A logic $\LL$ has a \textit{protodisjunction}\footnote{It is worth to remark that a more general notion of protodisjunction was introduced in  \cite[Convention 3.1]{CN13}, where protodisjunctions are understood as \textit{sets} of formulas in two variables.} if there is a binary term  $\lor$ such that $x \vdash_{\LL} x \lor  y$ and $y \vdash_{\LL} x \lor  y$.
\end{Definition}

Observe that every logic with theorems has a protodisjunction. This can be easily proved by showing that each theorem can be converted into a protodisjunction, by replacing each variable occurring in it by the variable $x$.
\begin{Theorem}\label{Thm:Chap2:Protodisjunction} If $\LL$ is Fregean and has a protodisjunction, then the parametrized equational translation $\btau(x, \vec{y}) = \{ x \lor y  \thickapprox x \}$ almost defines truth in $\ModstarL$.
\end{Theorem}
\begin{proof}
Thanks to Lemma \ref{Lem:Chap2:Transfer} it will be enough to prove that $\btau$ almost defines truth in $\LModstarL$. By Lemma \ref{Lem:Chap2:technical} this amounts to checking that
\[
\Gamma \vdash_{\LL}\varphi \text{ if and only if }\langle \varphi,  \varphi \lor \psi \rangle \in \leibniz \Gamma\text{ for every }\psi \in Fm
\]
for every $\Gamma \in \Th\LL \smallsetminus \{ \emptyset \}$ and $\varphi \in Fm$. For the ``only if'' part assume that $\Gamma \vdash_{\LL} \varphi$ and observe that $\varphi, \Gamma \sineq_{\LL}\Gamma, \varphi \lor \psi$, because $\lor$ is a protodisjunction. Since $\LL$ is Fregean, this implies that $\langle\varphi,  \varphi \lor \psi \rangle \in \frege_{\LL}\Gamma = \tarski_{\LL}\Gamma \subseteq \leibniz \Gamma$. For the ``if'' part choose $\psi \in \Gamma$. This can be done since $\Gamma \ne \emptyset$. Since $\lor$ is a protodisjunction, we have that $\Gamma \vdash_{\LL} \varphi \lor \psi$. By compatibility we conclude that $\Gamma \vdash_{\LL} \varphi$.
\end{proof}
Drawing consequences from this result, we obtain an essentially different proof of the following known result \cite[Theorem 14]{AFRM15}:
\begin{Corollary}\label{Cor:FregeanLogicsWithTheorems} A Fregean logic $\LL$ has theorems if and only if truth is  equationally definable in $\ModstarL$.\end{Corollary}
\begin{proof}
The ``if'' part was proven in Corollary \ref{Cor:Chap2:NoTheorem}. For the ``only if'' part, let $\LL$ be a Fregean logic with theorems. In particular, $\LL$ has a protodisjunction. Thus from Theorem \ref{Thm:Chap2:Protodisjunction} it follows that truth is almost parametrically equationally definable in $\ModstarL$. Together with Corollary \ref{Cor:Chap2:NoTheorem} and the fact that $\LL$ has theorems, this implies that truth is equationally definable in $\ModstarL$.
\end{proof}

Now we consider Fregean logics with a binary connective which behaves like a weak conjunction.
\begin{Definition} A logic $\LL$ has a \textit{protoconjunction} if there is a binary term $\land$ such that $\{x, y\} \vdash_{\LL} x \land y$, $\{x, x \land y\} \vdash_{\LL} y$ and $\{y, x \land y\} \vdash_{\LL} x$.
\end{Definition}

\begin{Theorem}\label{Thm:Chap2:Protoconjunction} If $\LL$ is Fregean and has a protoconjunction, then the parametrized equational translation $\btau(x, \vec{y}) = \{x \land y \thickapprox y \}$ almost defines truth in $\ModstarL$.
\end{Theorem}
\begin{proof}
Again it will be enough to prove that 
\[
\Gamma \vdash_{\LL}\varphi \text{ if and only if }\langle \psi,  \varphi \land \psi \rangle \in \leibniz \Gamma\text{ for every }\psi \in Fm
\]
for every $\Gamma \in \Th\LL \smallsetminus \{ \emptyset \}$ and $\varphi \in Fm$. For the ``only if'' part observe that if $\Gamma \vdash_{\LL} \varphi$, then $\psi, \Gamma \sineq_{\LL}\Gamma, \varphi \land \psi$ since $\land$ is a protoconjunction. Since $\LL$ is Fregean, this implies that $\langle \psi,  \varphi \land \psi \rangle \in \frege (\Th\LL)^{\Gamma} = \tarski_{\LL}\Gamma \subseteq \leibniz \Gamma$. For the ``if'' part choose $\psi \in \Gamma$. This can be done since $\Gamma \ne \emptyset$. By compatibility we have that $\Gamma \vdash_{\LL} \varphi \land \psi$. Since $\land$ is a protoconjunction, we conclude that $\Gamma \vdash_{\LL} \varphi$.
\end{proof}

\section{Implicit definability}

\begin{Definition} Truth is \textit{implicitly definable} in a class of matrices $\class{M}$ if the elements of $\class{M}$ are determined by their algebraic reducts, i.e., if $\langle \A, F\rangle, \langle \A, G\rangle\in \class{M}$, then $F = G$.
\end{Definition}

Similarly, we say that truth is \textit{almost} implicitly definable in a class of matrices $\class{M}$ when truth is implicitly definable in the class of non-almost trivial members of $\mathsf{M}$.

It is clear that (almost) implicit definability generalizes (almost) parametrized equational definability. Following our convention, we say that the \textit{truth sets} of a logic $\LL$ are (almost) \textit{implicitly definable}, when truth is (almost) implicitly definable in $\ModstarL$. The following observation is a straightforward application of the definitions:

\begin{Lemma}\label{Lem:Chap3:Implicit}
The truth sets of $\LL$ are (almost) implicitly definable if and only if $\leibniz^{\A} \colon \FFiLA \to \Con\A$ is (almost) injective for every algebra $\A$.
\end{Lemma}

When applied to classes of matrices, the concept of implicit definability has been inspired by Beth's definability theorem which states that in first-order logic explicit definability and implicit definability coincide. In particular, Blok and Pigozzi \cite{BP89}, Hermann \cite{He93a, He97}, Czelakowski and Jansana \cite{CzJa00} and Raftery \cite{Ra06a} proved a series of increasingly more general results, collectively called \textit{Beth's definability theorems}, whose main outcomes are summarized in the next result. Recall that a logic $\LL$ is \textit{protoalgebraic} \cite{BP86,Cz85,Cz86} if there is a set of formulas $\Delta(x, y)$ such that
\[
\emptyset \vdash_{\LL}\Delta(x, x) \quad x, \Delta(x, y) \vdash_{\LL} y.
\]
A logic is \textit{mono-unary} \cite{BlRe03,Ra06a}  when its language consists of just one symbol that, moreover, has arity $1$ (so that constants are not allowed). The following form of the theorem is due to Raftery  \cite[Corollary 29, Theorems 36 and 46]{Ra06a}.

\begin{Theorem}[Beth's definability]\label{Thm:Chap3:Beth}
The notions of implicit and equational definability coincide when applied to the class of Suszko-reduced models of a logic. Moreover, if a logic $\LL$ is either protoalgebraic or mono-unary and its truth predicates are implicitly definable, then they are equationally definable as well.
\end{Theorem}

Related to the study of Beth's definability theorems, another question inspired the research on the notion of implicit definability in matrix semantics. This is the so-called \textit{transfer problem} that asks whether the injectivity of the Leibniz operator transfers, in general, from the theories of a given logic to its deductive filters over arbitrary algebras. In other words the problem asks whether it is enough to check that truth is implicitly definable in $\LModstarL$, in order to conclude that the truth sets of $\LL$ are implicitly definable (cfr. Lemma \ref{Lem:Chap3:Implicit}). Czelakowski and Jansana  \cite{CzJa00} and Raftery \cite{Ra06a} proved respectively that the transfer problem has a positive answer both for protoalgebraic and mono-unary logics. In this contribution, we show that the transfer problem has a \textit{negative} solution in general solving \cite[Problem 1]{Ra06a}, see also \cite{DeMa05}. On the other hand, we prove that a positive answer can be recovered for logics expressed in a \textit{countable} language, generalizing Raftery's result on mono-unary logics. Moreover, we offer a new and more elementary proof of Jansana and Czelakowski's result on protoalgebraic logic. We begin by this last point:

\begin{Theorem}[Czelakowski and Jansana]\label{Thm:Chap3:ProtoalgebraicLogics}
For protoalgebraic\index{transfer of injectivity}\index{operator!Leibniz}\index{operator!(almost) injective}\index{logic!protoalgebraic}  logics the injectivity of the Leibniz operator transfers from theories to filters over arbitrary algebras.
\end{Theorem}
\begin{proof}[Alternative proof]
Suppose that $\LL$ is protoalgebraic and that $\leibniz$ is injective over $\Th\LL$. We will show that $\leibniz$ is also completely order-reflecting over $\Th\LL$. Consider a family $\Sigma \cup \{ \Gamma \} \subseteq \Th\LL$ such that $\bigcap \{ \leibniz \Gamma' : \Gamma' \in \Sigma \} \subseteq \leibniz \Gamma$. It is well known that the fact that a logic $\LL'$ is protoalgebraic is equivalent to the fact that $\leibniz^{\A}F = \tarski_{\LL'}F$ for every algebra $\A$ and filter $F \in \FFi_{\LL'}A$ \cite[Theorem 6.7]{AAL-AIT-f}. Keeping this in mind, we obtain that:
\[
\leibniz \bigcap \Sigma = \tarski_{\LL} \bigcap \Sigma \subseteq \bigcap \{ \leibniz \Gamma' : \Gamma' \in \Sigma \} \subseteq \leibniz \Gamma
\]
Moreover, a logic $\LL'$ is protoalgebraic if and only if $\leibniz^{\A} \colon \FFi_{\LL'}\A \to \Con\A$ commutes with arbitrary meets for every algebra $\A$ \cite[Proposition 6.14]{AAL-AIT-f}. Applying this property to $\LL$, we obtain that
\[
\leibniz \bigcap \Sigma  = \leibniz \Gamma \cap \leibniz \bigcap \Sigma = \leibniz (\Gamma \cap\bigcap \Sigma).
\]
Since $\leibniz$ is injective over $\Th\LL$, we conclude that $\bigcap \Sigma = \Gamma \cap\bigcap \Sigma$, that is $\bigcap \Sigma \subseteq \Gamma$. With an application of Theorem \ref{Thm:Chap2:LeibnizEquational} we are done.
\end{proof}

In fact the proof presented above shows more than what is contained in the statement of Theorem \ref{Thm:Chap3:ProtoalgebraicLogics}. Namely, it establishes the part of the Beth's Definability Theorem \ref{Thm:Chap3:Beth} concerning protoalgebraic logics. Now we turn to prove that the transfer problem has a positive solution for logics expressed in a countable language. To this end, we will make use of the following technical result which generalizes \cite[Proposition 0.7.6]{Cz01}.

\begin{Lemma}\label{Lem:Chap3:Linguistic} Let\index{matrix!reduced} $\kappa$ be an infinite cardinal larger or equal to the cardinality of the language $\vert \mathscr{L} \vert$ and let $\langle \A, F \rangle$ and $\langle \A, G \rangle$ be a pair of reduced matrices. Every $\kappa$-generated subalgebra $\C$ of $\A$ can be extended to another $\kappa$-generated subalgebra $\B$ of $\A$ such that the matrices $\langle \B, F \cap B \rangle$ and $\langle \B, G \cap B \rangle$ are reduced.
\end{Lemma}
\begin{proof}
First observe that $\kappa$-generated algebras are of cardinality $\leq \kappa$, since $\kappa$ is infinite and larger or equal to $\vert \mathscr{L} \vert$. Then we define recursively an infinite family of subsets of $A$. We begin with $X_{0} \coloneqq C$. To define $X_{n+1}$, we go through the following construction: for every pair of different elements $a, b\in \textup{Sg}^{\A}(X_{n})$ we pick two finite sequences $\vec{c}$ and $\vec{d}$ of elements of $A$ for which there is a pair of formulas $\varphi(x, \vec{y})$ and $\psi(x, \vec{z})$ such that
\begin{align*}
\varphi^{\A}(a, \vec{c})\in F &\Longleftrightarrow \varphi^{\A}(b, \vec{c})\notin F\\
\psi^{\A}(a, \vec{d})\in G &\Longleftrightarrow \psi^{\A}(b, \vec{d})\notin G.
\end{align*}
The existence of the sequences $\vec{c}$ and $\vec{d}$ is ensured by point 1 of Lemma \ref{Lem:Polynomial} together with the fact that $\langle a, b \rangle \notin \textup{Id}_{\A} = \leibniz^{\A}F = \leibniz^{\A}G$.  We then let $Y_{n}$ be the set of all elements in the sequences constructed in this way. Finally we set
\begin{align*}
X_{n+1} \coloneqq X_{n} \cup Y_{n}.
\end{align*}

Now consider the union $\bigcup_{n\in \omega } \textup{Sg}^{\A}(X_{n})$. It is easy to prove that it is the universe 
of a subalgebra $\B$ of $\A$. Clearly $\B$ extends $\C$, since $X_{0} = C$. We claim that $\langle \B, F \cap B\rangle$ is reduced. To prove this, consider two different $a, b\in B$. There is $n \in \omega$ such that $a, b\in \textup{Sg}^{\A}(X_{n})$. By definition of $X_{n+1}$, we know that there is a finite sequence 
$\vec{c}$ of elements of $X_{n+1}$ and a formula $\varphi(x, \vec{y})$ such 
that $\varphi^{\A}(a, \vec{c})\in F \Longleftrightarrow \varphi^{\A}(b, \vec{c}) \notin F$. Since $a, b, \vec{c} \in B$ and $\B$ is a subalgebra of $\A$, 
we conclude that $\varphi^{\B}(a, \vec{c}), \varphi^{\B}(b, \vec{c}) \in B$ and, finally, that
\begin{align*}
\varphi^{\B}(a, \vec{c}) \in F \cap B \Longleftrightarrow \varphi^{\B}(b, \vec{c}) \notin F \cap B. 
\end{align*}
By point 1 of Lemma \ref{Lem:Polynomial} we conclude that $\langle a, b \rangle \notin \leibniz^{\B}(F \cap B)$ and therefore that $\leibniz^{\B}(F \cap B) = 0_{\B}$. 
This concludes the proof of our claim. An analogous argument yields that the the matrix $\langle \B, G \cap B \rangle$ is reduced too.

It only remains to show that $\B$ is $\kappa$-generated. We begin by showing inductively that $\vert X_{n} \vert \leq \kappa$ 
for every $n\in \omega$. For $n=0$ we have that $X_{0}= C$. Recall from the assumption that $\C$ 
is $\kappa$-generated and, therefore, of cardinality $\leq \kappa$. For the $n+1$ case observe that, by the inductive hypothesis, $\vert X_{n} \vert \leq \kappa$. Then $\textup{Sg}^{\A}(X_{n})$ is $\kappa$-generated and again of cardinality $\leq \kappa$. Now observe 
that, while constructing $Y_{n}$, we added to $X_{n}$ at most a finite number of elements for every pair of 
elements of $\textup{Sg}^{\A}(X_{n})$. Therefore the cardinality of $Y_{n}$ can be bounded above by 
$\aleph_{0} \cdot \vert \textup{Sg}^{\A}(X_{n}) \vert  \cdot \vert \textup{Sg}^{\A}(X_{n}) \vert$. In particular, 
this yields that
\begin{align*}
\vert Y_{n} \vert \leq \aleph_{0} \cdot \vert \textup{Sg}^{\A}(X_{n}) \vert  \cdot \vert \textup{Sg}^{\A}(X_{n}) \vert  \leq \aleph_{0} \cdot (\kappa \cdot \kappa) \leq \kappa.
\end{align*}
Since $X_{n+1} = X_{n} \cup Y_{n}$ is the union of two sets of cardinality smaller or equal to $\kappa$, we conclude that $\vert X_{n+1}\vert \leq \kappa$. This concludes our proof by 
induction. Then let $n\in \omega$. The fact that $\vert X_{n} \vert \leq \kappa$ implies that $\textup{Sg}^{\A}(X_{n})$is  of cardinality $\leq \kappa$. Thus $B = \bigcup_{n\in \omega} \textup{Sg}^{\A}(X_{n})$ is the union of countably many sets of cardinality smaller of equal to $\kappa$. Since $\kappa$ is infinite, we conclude that $\B$ has cardinality  $\leq \kappa$, hence a fortiori it is $\kappa$-generated.
\end{proof}

We are now ready to prove our transfer result:

\begin{Theorem}\label{Thm:Chap3:TransferCountable}
For logics expressed\index{language!countable}\index{transfer of injectivity}\index{operator!Leibniz}\index{operator!(almost) injective}  in a countable language the (almost) injectivity of the Leibniz operator transfers from theories to filters over arbitrary algebras.
\end{Theorem}

\begin{proof}
We apply Lemma \ref{Lem:Chap3:Implicit}. Consider two reduced models $\langle \A, F \rangle$ and $\langle \A, G \rangle$ of $\LL$. We have to prove that $F = G$. By symmetry it is enough to prove that $F \subseteq G$. Consider an element $a \in F$ and let $\C$ be the subalgebra of $\A$ generated by $\{ a \}$. Clearly $\C$ is countably generated. Since $\vert \mathscr{L} \vert \leq \aleph_{0}$, we can apply Lemma \ref{Lem:Chap3:Linguistic} and extend $\C$ to a countably generated subalgebra $\B$ of $\A$ such that both $\langle \B, F \cap B \rangle$ and $\langle \B, G \cap B \rangle$ are reduced matrices. Clearly they are both models of $\LL$. Since $\B$ is countably generated, we can choose a surjective homomorphism $h \colon \Fm \to \B$. Then we define $\Gamma \coloneqq h^{-1}[F \cap B]$ and $\Gamma' \coloneqq h^{-1}[G \cap B]$. We have $\Gamma, \Gamma \in \Th\LL$. With an application of Lemma \ref{Lem:Chap1:Inverse} we obtain
\begin{eqnarray*}
\leibniz \Gamma &=& \leibniz h^{-1}[F \cap B] = h^{-1}\leibniz^{\B}(F \cap B ) = h^{-1} 0_{\B}\\
&=& h^{-1}\leibniz^{\B}(G \cap B ) = \leibniz h^{-1}[G \cap B] = \leibniz \Gamma'.
\end{eqnarray*}
Together with the assumption, this implies that $\Gamma = \Gamma'$. In particular, this implies that $a \in F \cap B = G \cap B \subseteq G$. This concludes the proof that $F \subseteq G$ and therefore we are done. The \textit{almost} case follows by restricting the proof to non-empty filters.
\end{proof}

\begin{Corollary}[Raftery]
For mono-unary logics the injectivity of the Leibniz operator transfers from theories to filters over arbitrary algebras.
\end{Corollary}

The following problem was raised by J. Gil-F\'erez in a private communication:

\begin{problem}
Is it possible to view Theorems \ref{Thm:Chap3:ProtoalgebraicLogics} and \ref{Thm:Chap3:TransferCountable} as special instances of some more general phenomenon?
\end{problem}

\section{Failure of the transfer of injectivity}\label{Sec:Failure}

As we mentioned, if we move our attention to logics expressed in uncountable languages, it is possible to construct examples where the injectivity of the Leibniz operator does not transfer from theories to deductive filters over arbitrary algebras. We devote this section to the description of such an example.

Let $\oper{R}$ be the set\index{transfer of injectivity}\index{operator!Leibniz}\index{operator!(almost) injective}\index{filters!deductive}
\index{language!uncountable} of real numbers. We consider the algebraic type that consists of a set of binary connectives $\{ \multimap_{i} : i \in \oper{R} \smallsetminus \{ 1, 2 \} \}$, a set of constants $\{ c_{i} : i \in \oper{R} \smallsetminus \{ 2 \} \}$ and a unary connective $\Box$. Then let $\A$ be the algebra with universe
\begin{align*}
A \coloneqq ((\oper{R} \smallsetminus \{ 2 \}) \times \{ 0 \}) \cup (\oper{R} \times \{ 1 \})
\end{align*}
and operations defined as follows:
\begin{eqnarray*}
\langle a_{1}, a_{2} \rangle \multimap_{i} \langle b_{1}, b_{2}\rangle &\coloneqq&
\begin{cases}
 \langle 1, 1 \rangle & \textrm{if $a_{1} = b_{1} = i$, $a_{2} = 0$ and $b_{2}=1$}\\
 \langle 1, 0 \rangle & \textrm{otherwise}
\end{cases}\\  
\Box \langle a_{1}, a_{2} \rangle &\coloneqq&
\begin{cases}
 \langle 1, 1 \rangle & \textrm{if either $a_{2}=0$ or ($a_{2}= 1$ and $a_{1} < 2$)}\\
 \langle 1, 0 \rangle & \textrm{otherwise}
\end{cases}\\
c_{i} &\coloneqq& \langle i, 1 \rangle
\end{eqnarray*}
for every $\langle a_{1}, a_{2} \rangle, \langle b_{1}, b_{2}\rangle \in A$. To simplify the notation put
\[
F \coloneqq   \oper{R} \times \{ 1 \} \text{ and } G \coloneqq (\oper{R} \times \{ 1 \}) \smallsetminus\{ \langle 2, 1 \rangle \}.
\]
We consider the logic $\LL$ determined by the pair of matrices $\langle \A, F\rangle$ and $\langle \A, G\rangle.$

\begin{fact}\label{Fact:NotImplicitlyDefinableInModstar} $\leibniz^{\A}$ is not injective over $\FFiL\A$.
\end{fact}
\begin{proof}
We claim that for every pair of different $\langle a_{1}, a_{2} \rangle, \langle b_{1}, b_{2}\rangle \in A$, there is a polynomial function $p(z)$ of $\A$ that satisfies one of the following conditions: 
\begin{align*}
p\langle a_{1}, a_{2} \rangle \in F \cap G &\text{ and } p \langle b_{1}, b_{2}\rangle \notin F \cup G\\
p\langle b_{1}, b_{2} \rangle \in F \cap G &\text{ and } p \langle a_{1}, a_{2}\rangle \notin F \cup G.
\end{align*}
\noindent We split the proof of the claim in three main cases: 

\benormal
\item $a_{2} \ne b_{2}$.
\item $a_{2} = b_{2} = 0$.
\item $a_{2} = b_{2} = 1$.
\enormal
1. Assume w.l.o.g. that $a_{2}= 0$ and $b_{2}= 1$. If $b_{1} \ne 2$, then $\langle a_{1}, a_{2} \rangle \notin F \cup G$ and $\langle b_{1}, b_{2}\rangle \in F \cap G$. Then suppose that $b_{1} = 2$. We have that
\begin{align*}
\Box \langle a_{1}, a_{2}\rangle = \langle 1, 1 \rangle \in F \cap G \text{ and }\Box \langle b_{1}, b_{2}\rangle = \langle 1, 0\rangle \notin F \cup G.
\end{align*}

2. Since $\langle a_{1}, a_{2}\rangle \ne \langle b_{1}, b_{2}\rangle$, we have that either $a_{1} \ne 1$ or $b_{1} \ne 1$. Assume w.l.o.g. that $a_{1} \ne 1$. We consider the unary polynomial function
\begin{align*}
p(z) \coloneqq z \multimap_{a_{1}}\langle a_{1}, 1 \rangle.
\end{align*}
Observe that the operation $\multimap_{a_{1}}$ exists, since $a_{1} \notin \{ 1, 2 \}$. Clearly we have that $p\langle a_{1}, a_{2}\rangle \in F \cap G$, while $p\langle b_{1}, b_{2}\rangle \notin F \cup G$. 

3. We have two subcases: either $\{ a_{1}, b_{1} \}\ne \{ 1, 2 \}$ or $\{ a_{1}, b_{1} \} = \{ 1, 2 \}$. Suppose that  $\{ a_{1}, b_{1} \}\ne \{ 1, 2 \}$. Since $\langle a_{1}, a_{2}\rangle \ne \langle b_{1}, b_{2}\rangle$ and $a_{2} = b_{2}$, we know that $a_{1} \ne b_{2}$. Together with the fact that $\{ a_{1}, b_{1} \}\ne \{ 1, 2 \}$, this implies that either $a_{1} \notin \{ 1, 2 \}$ or $b_{1} \notin \{ 1, 2 \}$. Assume w.l.o.g. that $a_{1} \notin \{ 1, 2 \}$. Then we can safely consider the polynomial function
\begin{align*}
p(z) \coloneqq \langle a_{1}, 0\rangle \multimap_{a_{1}} z.
\end{align*}
It is easy to see that $p \langle a_{1}, a_{2}\rangle \in F \cap G$, while $p\langle b_{1}, b_{2}\rangle \notin F \cup G$. Then we consider the case where $\{ a_{1}, b_{1} \} = \{ 1, 2 \}$. Assume w.l.o.g. that $a_{1} = 1$ and $b_{1} =2$. We have that
\begin{align*}
\Box \langle a_{1}, a_{2}\rangle = \langle 1, 1\rangle\in F \cap G\text{ and }\Box \langle b_{1}, b_{2}\rangle = \langle 1, 0 \rangle \notin G \cup F.
\end{align*}
This establishes our claim. From Lemma \ref{Lem:Polynomial} it follows that $\leibniz^{\A}F = \leibniz^{\A}G$. Since $F \ne G$ we are done.
\end{proof}

\begin{fact}\label{Fact:Chap3:Multimap} Consider $\Gamma \in \Th\LL$, $\varphi \in \Gamma$ and a formula of the form $\alpha(z) \multimap_{i} \beta(z)$ in which $z$ actually occurs. If $\Gamma \vdash_{\LL}\alpha(\varphi) \multimap_{i} \beta(\varphi)$, then $\langle \varphi, c_{i}\rangle \in \leibniz \Gamma$.
\end{fact}

\begin{proof}
Suppose that $\Gamma \vdash_{\LL}\alpha(\varphi) \multimap_{i} \beta(\varphi)$. By Lemma \ref{Lem:Polynomial} it will be enough to prove that $h(\varphi) = h(c_{i})$ for every homomorphism $h \colon \Fm \to \A$ such that $h[\Gamma] \subseteq \oper{R} \times \{ 1 \}$. Then consider an homomorphism $h$ of this kind. Since $h[\Gamma] \subseteq \oper{R} \times \{ 1 \}$ we have that $h(\alpha(\varphi) \multimap_{i} \beta(\varphi)) \in \oper{R} \times \{ 1 \}$. Looking at the definition of $\multimap_{i}$, it is easy to see that this happens only if $h(\alpha(\varphi) \multimap_{i} \beta(\varphi)) = \langle 1, 1 \rangle$. In particular, this is to say that
\[
h\alpha(\varphi) = \langle i, 0\rangle\text{ and }h\beta(\varphi) = \langle i, 1\rangle.
\]
Looking at the definition of the basic operations of $\A$ and keeping in mind that $i \in \oper{R}\smallsetminus \{ 1, 2 \}$, it is possible to see that $\alpha(\varphi)$ must be a variable and that $\beta(\varphi)$ must be either a variable or $c_{i}$. Then we have cases:
\begin{equation}\label{Eq:Chap3:Multimap}
\text{either }\alpha(\varphi) \multimap_{i} \beta(\varphi) = x \multimap_{i} y \text{ or } \alpha(\varphi) \multimap_{i} \beta(\varphi) = x \multimap_{i} c_{i}
\end{equation}
for some variables $x$ and $y$. Now, from the assumption we know that $z$ occurs in $\alpha(z) \multimap_{i} \beta(z)$. We claim that $z$ does not appear really in $\alpha(z)$. To prove this, suppose the contrary towards a contradiction. By (\ref{Eq:Chap3:Multimap}) we would have that $\varphi = x$. Then
\[
h\varphi = h \alpha(\varphi) = \langle i, 0 \rangle \notin \oper{R} \times \{ 1 \}.
\]
But this contradicts the fact that $\Gamma \vdash_{\LL} \varphi$, establishing the claim. In particular, the claim implies that $z$ occurs in $\beta(z)$. Together with (\ref{Eq:Chap3:Multimap}), this means that $\beta(\varphi) = \varphi$. This easily implies that
\[
h(\varphi) = h\beta(\varphi) = \langle i, 1 \rangle = h(c_{i})
\]
establishing the fact.
\end{proof}

\begin{fact}\label{Fact:Chap3:VariablesConstants}
If $\varphi$ is neither a variable nor a constant, then $\emptyset \vdash_{\LL}\Box \varphi$.
\end{fact}
\begin{proof}
Observe that $\Box \Box x$ and $\Box( x \multimap_{i} y)$ are theorems of $\LL$, for every $i \notin \{ 1, 2 \}$. This follows directly from the definition of $\LL$. In particular, this implies that also $\Box \varphi$ is a theorem.
\end{proof}

\begin{fact}\label{Fact:Chap3:BoxEqual} For every $\Gamma \cup \{ \varphi \} \subseteq Fm$ one of the following conditions hold:
\benormal
\item For every $i < 2$ and formula $\alpha(z)$ in which $z$ actually occurs: $\Gamma \vdash_{\LL}\Box\alpha(\varphi) \Longleftrightarrow \Gamma \vdash_{\LL}\Box\alpha(c_{i})$.
\item For every $i > 2$ and formula $\alpha(z)$ in which $z$ actually occurs: $\Gamma \vdash_{\LL}\Box\alpha(\varphi) \Longleftrightarrow \Gamma \vdash_{\LL}\Box\alpha(c_{i})$.
\enormal
\end{fact}

\begin{proof}
First consider the case in which $\Gamma \vdash_{\LL} \Box \alpha(\varphi)$ for every formula $\alpha(z)$ in which $z$ really occurs. We want to prove that condition 1 is satisfied. Then consider $c_{i}$ with $i < 2$ and a formula $\alpha(z)$ in which $z$ really occurs. If $\alpha(z)$ is not a variable, then $\Box \alpha(c_{i})$ is a theorem by Fact \ref{Fact:Chap3:VariablesConstants}. Then consider the case where $\alpha(z)$ is a variable. Clearly $\alpha = z$. Also in this case $\Box \alpha(c_{i})$ is a theorem. Thus we conclude that condition 1 is satisfied.

Then consider the case where there is at least one formula $\alpha(z)$, in which $z$ really occurs, such that $\Gamma \nvdash_{\LL} \Box \alpha(\varphi)$. From Fact \ref{Fact:Chap3:VariablesConstants} it follows that $\alpha(z) = z$ and, therefore, that $\Gamma \nvdash_{\LL}\varphi$. Together with Fact \ref{Fact:Chap3:VariablesConstants} this implies that for every formula $\beta(z)$, in which $z$ really occurs, we have that
\[
\Gamma \vdash_{\LL}\Box \beta(\varphi) \Longleftrightarrow \beta \ne z.
\]
Then consider a constant $c_{i}$ with $i >2$ and a formula $\beta(z)$ in which $z$ really occurs. From Fact \ref{Fact:Chap3:VariablesConstants} and from the definition of $\Box$ it follows that $\Gamma \vdash_{\LL}\Box \beta(c_{i})$ if and only if $\beta \ne z$. Thus condition 2 is satisfied.
\end{proof}

\begin{fact}\label{Fact:Chap3:ThereIsAReal} For every $\Gamma \in \Th\LL$ and $\varphi \in \Gamma$, there is $i \in \oper{R} \smallsetminus \{ 2 \}$ such that $\langle \varphi, c_{i}\rangle \in \leibniz \Gamma$.
\end{fact}
\begin{proof}
If $\Gamma$ is the inconsistent theory, then $\leibniz \Gamma = Fm \times Fm$ and, therefore, we are done. Then consider the case where $\Gamma$ is consistent. Suppose towards a contradiction that $\langle \varphi, c_{i}\rangle \notin \leibniz \Gamma$ for every $i \in \oper{R} \smallsetminus \{ 2 \}$. From Lemma \ref{Lem:Polynomial} it follows that for every $i \in \oper{R} \smallsetminus \{ 2 \}$ there is a formula $p(z)$ such that
\begin{equation}\label{Eq:ThereIsAPolynomial}
\Gamma \vdash_{\LL}p(\varphi) \Longleftrightarrow \Gamma \nvdash_{\LL}p(c_{i}).
\end{equation}
By Fact \ref{Fact:Chap3:BoxEqual} one of the following conditions hold:
\benormal
\item For every $i < 2$ and formula $\alpha(z)$ in which $z$ actually occurs: $\Gamma \vdash_{\LL}\Box\alpha(\varphi) \Longleftrightarrow \Gamma \vdash_{\LL}\Box\alpha(c_{i})$.
\item For every $i > 2$ and formula $\alpha(z)$ in which $z$ actually occurs: $\Gamma \vdash_{\LL}\Box\alpha(\varphi) \Longleftrightarrow \Gamma \vdash_{\LL}\Box\alpha(c_{i})$.
\enormal
Assume that condition 1 holds (the proof for case 2 is analogous). Then consider $1 \ne i < 2$. 
There is a polynomial function $p(z)$ that satisfies (\ref{Eq:ThereIsAPolynomial}). Thus $z$ actually occurs in $p(z)$. By condition 1, we know that the main connective of $p(z)$ cannot be $\Box$. Therefore $p(z) = \alpha(z) \multimap_{j} \beta(z)$ for some $j \in \oper{R} \smallsetminus \{ 1 , 2 \}$ and formulas $\alpha$ and $\beta$. Together with Fact \ref{Fact:Chap3:Multimap} and $\langle \varphi, c_{i}\rangle \notin \leibniz \Gamma$, this implies that
\[
\Gamma \vdash_{\LL}\alpha(c_{i}) \multimap_{j} \beta(c_{i}).
\]
Since $\Gamma$ is consistent, there is a homomorphism $h \colon \Fm \to \A$ such that $h[\Gamma] \subseteq \oper{R} \times \{ 1 \}$. We have that $h(\alpha(c_{i}) \multimap_{j} \beta(c_{i})) \in \oper{R} \times \{ 1 \}$. But this is to say that
\[
h \alpha(c_{i}) = \langle j, 0\rangle \text{ and }h \beta(c_{i}) = \langle j, 1\rangle.
\]
From the definition of $\multimap_{j}$ it follows that $\alpha(\varphi)$ must be a variable $y_{i}$ and that $\beta(\varphi)$ must be either a variable or $c_{j}$ (this is because $i \ne 1$). Since $z$ actually occurs in $\alpha(z) \multimap_{j} \beta(z)$ we conclude that $\alpha(c_{i}) \multimap_{j} \beta(c_{i}) = y_{i} \multimap_{j} c_{i}$. Keeping in mind that $h(y_{i}\multimap_{j} c_{i}) \in \oper{R} \times \{ 1 \}$, we obtain that $j = i$. Therefore we conclude that $\Gamma \vdash_{\LL} y_{i} \multimap_{i} c_{i}$.  Now, we proved that for every $1 \ne i < 2$ there is a variable $y_{i}$ such that
\begin{equation}\label{Eq:Chap3:FinalStep}
\Gamma \vdash_{\LL} y_{i} \multimap_{i} c_{i}.
\end{equation}
Consider the homomorphism $h$ above. From (\ref{Eq:Chap3:FinalStep}) it follows that $h (y_{i}) = \langle i, 0\rangle$ for every $i < 2$. But this contradicts the fact that there are uncountably many reals smaller than $2$ and different from $1$ and only countably many variables in $Fm$.
\end{proof}

\begin{fact}\label{Fact:NotTrasfer} $\leibniz$ is order-reflecting (and therefore injective) over $\Th\LL$.
\end{fact}

\begin{proof}
Consider two theories $\Gamma, \Gamma' \in \Th\LL$ such that $\leibniz \Gamma \subseteq \leibniz \Gamma'$. Then pick $\varphi \in \Gamma$ (we can always do this, since $\LL$ has theorems). By Fact \ref{Fact:Chap3:ThereIsAReal} there is a constant $c_{i}$ such that $\langle \varphi, c_{i}\rangle \in \leibniz \Gamma \subseteq \leibniz \Gamma'$. Since $\emptyset \vdash_{\LL} c_{i}$, by compatibility we obtain that $\varphi \in \Gamma'$. Hence $\leibniz \Gamma \subseteq \leibniz \Gamma'$ as desired.
\end{proof}

The outcome of the work done in this section can be summarized as follows:
\begin{Theorem}
The injectivity of the Leibniz operator does not transfer in general from theories to filters over arbitrary algebras.
\end{Theorem}

\section{Small truth sets}\label{Sec:SmallTruthSets}

There is at least another definability condition that fits into the framework of the Leibniz hierarchy, and that until now was not recognized in the literature.

\begin{Definition}
Let $\class{M}$ be a class\index{truth!(almost) small} of matrices and $\LL$ the logic it defines. Truth is \textit{small} in $\class{M}$ if the truth set $F$ is the smallest non-empty deductive filter of $\LL$ over $\A$, for every $\langle \A, F \rangle \in \mathsf{M}$.
\end{Definition}

As usual, we say that the \textit{truth sets } of logic $\LL$ are \textit{small}, if truth is small in $\ModstarL$. We have the following:
\begin{Lemma}\label{Lem:BasicWisdom} \
\benormal
\item If the truth sets of a logic are (almost) small, then they are implicitly definable as well.
\item If the truth sets of a logic are (almost parametrically equationally) equationally definable, then they are (almost) small as well.
\enormal
\end{Lemma}
\begin{proof}
1. Suppose that truth is almost small in $\ModstarL$. Then consider two non-almost trivial reduced models $\langle \A, F \rangle$ and $\langle \A, G\rangle$ of $\LL$. Since truth is almost small in $\ModstarL$, we have that both $F$ and $G$ are the smallest non-empty deductive filter of $\LL$ over $\A$. In particular, this means that $F = G$.

2. Suppose that truth is almost parametrically equationally definable in $\ModstarL$. Then consider a non-almost trivial reduced model $\langle \A, F \rangle$ of $\LL$ and any $G \in \FFiL\A \smallsetminus \{ \emptyset \}$. We have that $\leibniz^{\A}F = \textup{Id}_{\A} \subseteq \leibniz^{\A}G$. With an application of Theorem \ref{Thm:Chap2:LeibnizUniversal} we obtain that $F \subseteq G$. Since $F \ne \emptyset$, we conclude that $F$ is the smallest non-empty deductive filter of $\LL$ over $\A$.
\end{proof}

\label{Comments} The fact that the notion of smallness is strictly stronger than the one of implicit definability, when referred the truth sets of a logic, was first proved by the \cite[Example 2]{Ra06a} of Raftery.\footnote{Observe that in the paper \cite{Ra06a} the notion of \textit{smallness} is not considered on his own right, so that the reader will not find there an explicit statement separating the concept of smallness from that of implicit definability. What is proved in the quoted \cite[Example 2]{Ra06a} is the existence of a logic $\LL$ such that $\leibniz^{\A} \colon \FFiLA \to \Con\A$ is injective, but not necessarily order-reflecting, for every algebra $\A$. In the light of Theorem \ref{Lem:Chap3:Reflecting}, this is enough to separate the concept of implicit definability from that of smallness.} Another example of logic whose truth sets are implicitly definable, but not small, is the $\langle \Box, 1\rangle$-fragment of the local consequence of the modal system $\mathcal{S}4$ \cite[Example 3.4]{TM16phd}. Moreover, the fact that the notion of equational definability is strictly stronger than the one of smallness, when referred to the truth sets of a logic, is proved in Example \ref{Ex:Chap3:Small2}.

The next result shows that the class of logics whose truth sets are small belongs to the Leibniz hierarchy, in the sense that it is characterized by a property of the Leibniz operator:

\begin{Theorem}\label{Lem:Chap3:Reflecting}
Truth is (almost) small in $\ModstarL$ if and\index{operator!Leibniz}\index{operator!(almost) order reflecting}\index{filters!deductive}\index{truth!(almost) small}\index{model!reduced} only if $\leibniz^{\A}$ is (almost) order reflecting over $\FFiLA$ for every algebra $\A$.
\end{Theorem}
\begin{proof}
We begin by the ``only if'' part. Consider $F, G \in \FFiLA \smallsetminus \{ \emptyset 	\}$ such that $\leibniz^{\A}F \subseteq \leibniz^{\A}G$. Let $h \colon \A/\leibniz^{\A} F \to \A/\leibniz^{\A} G$ be the natural surjection. We have that $h^{-1}[G/\leibniz^{\A}G] \in \FFiL(\A/\leibniz^{\A}F)$. By the assumption we know that $F/ \leibniz^{\A}F$ is the smallest non-empty deductive filter of $\LL$ on $\A / \leibniz^{\A}F$. Since $h^{-1}[G/\leibniz^{\A}G] \ne \emptyset$, we conclude that $F/ \leibniz^{\A}F \subseteq h^{-1}[G/\leibniz^{\A}G]$. Then let $a \in F$. We have that $a/ \leibniz^{\A}F\subseteq h^{-1}[G/\leibniz^{\A}G]$ and, therefore, $a \in G$.

Then we turn to check the ``if'' part. Consider a non-almost trivial reduced model $\langle \A, F\rangle$ of $\LL$ and a filter $G \in \FFiLA \smallsetminus \{ \emptyset \}$. We have that $\leibniz^{\A}F = 0_{\A} \subseteq \leibniz^{\A}G$. By the assumption we obtain $F \subseteq G$. Since $F \ne \emptyset$, we conclude that $F$ is the smallest non-empty deductive filter of $\LL$ on $\A$.
\end{proof}

Related to the study of the transfer problem, it is natural to ask whether the order-reflection of the Leibniz operator transfers from theories to filters over arbitrary algebras. This is not the case in general: the logic $\LL$ described in Section \ref{Sec:Failure} is such that $\leibniz \colon \Th\LL \to \Con\Fm$ is order-reflecting (Fact \ref{Fact:NotTrasfer}), while $\leibniz^{\A} \colon \FFiLA \to \Con\A$ is not injective (and, therefore, not order-reflecting) in general. Nevertheless, a natural adaptation of the proof of Theorem  \ref{Thm:Chap3:TransferCountable} yields the following result:

\begin{Theorem}\label{Thm:Chap3:TransferReflecting}
For logics expressed in a countable language the (almost)\index{language!countable}\index{operator!Leibniz}\index{operator!(almost) order reflecting}\index{filters!deductive}\index{truth!(almost) small}\index{model!reduced}  order-reflection of the Leibniz operator transfers from theories to filters over arbitrary algebras.
\end{Theorem}
\begin{proof}
We apply Theorem \ref{Lem:Chap3:Reflecting}. Suppose that $\leibniz$ is order-reflecting over $\Th\LL$. This easily implies that $\LL$ has theorems and, therefore, that its deductive filters are non-empty. Consider a reduced model $\langle \A, F \rangle$ of $\LL$. We have to prove that $F$ is the smallest non-empty deductive filter of $\LL$ over $\A$. Then consider $G \in \FFiLA \smallsetminus \{ \emptyset \}$ and $a\in F$. We know that $G \ne \emptyset$. Thus we can choose an element $b\in G$. We let $\C$ be the subalgebra of $\A$, generated by $\{a, b\}$. Since $\C$ is finitely generated, we can apply Lemma \ref{Lem:Chap3:Linguistic} and extend it to a countably generated subalgebra $\B$ of $\A$ such that $\langle \B, F \cap B \rangle$ is a reduced model of $\LL$. 

Since $\B$ is countable, there is a surjective homomorphism $h \colon \Fm \to \B$. Then let $\Gamma \coloneqq h^{-1}[F \cap B]$ and $\Gamma' \coloneqq h^{-1}[G \cap B]$. Notice that $F\cap B$ and $G\cap B$ are non-empty, and hence $\Gamma$ and $\Gamma'$ are non-empty as well. Clearly $\Gamma, \Gamma' \in \Th\LL$. From Lemma \ref{Lem:Chap1:Inverse} we obtain that
\begin{align*}
\leibniz\Gamma = \leibniz h^{-1}[F \cap B] = h^{-1}\leibniz^{\B}( F \cap B )= h^{-1}0_{\B} = \textup{Ker}(h).
\end{align*}
Moreover, since $\textup{Ker}(h)$ is compatible with $\Gamma'$, we know that $\leibniz\Gamma =\textup{Ker}(h) \subseteq \leibniz\Gamma'$. Hence we can apply the assumption and conclude that $\Gamma \subseteq \Gamma'$. Together with the fact that $h$ is surjective, this implies that
\begin{align*}
F \cap B = hh^{-1}[ F\cap B] = h[\Gamma] \subseteq h[\Gamma'] = hh^{-1}[G \cap B] \subseteq G \cap B.
\end{align*}
Therefore we obtain that $a \in F \cap B \subseteq G \cap B \subseteq G$. This shows that $F \subseteq G$.
\end{proof}

\begin{exa}[\textsf{Small truth sets}]\label{Ex:Chap3:Small2}
We describe a logic $\LL$, whose truth sets are small but not equationally definable. To this end, let $\LL$ be the logic, expressed in the language $\langle \Box, 1 \rangle$ of type $\langle 1, 0 \rangle$, axiomatized by the following Hilbert-style rules:
\begin{align*}
\emptyset \vdash 1 \quad \emptyset \vdash \Box 1 \quad \Box \Box x \vdash y.
\end{align*}
We will prove that truth is small, but not equationally definable, in $\ModstarL$. To this end, let $\A_{4} = \langle \{ a, b, c,  1 \}, \Box, 1\rangle$ be the algebra where $\Box$ is defined for every $x \in A_{4}$ as follows:
\begin{align*}
\Box p = \left\{ \begin{array}{ll}
a & \textrm{if $p \in \{ 1, c\}$}\\
 b & \textrm{otherwise.}\\
\end{array} \right.
\end{align*}
Let $\A_{3}$ be the subalgebra of $\A_{4}$ with universe $\{1, a, b\}$. 
\begin{fact}\label{Fact:CUIso}
$\ModstarL$ is the closure under isomorphism of $\langle \A_{4}, \{ 1, a \}\rangle$, $\langle \A_{3}, \{1, a \}\rangle$ and $\langle {\bf 1}, \{ 1 \} \rangle$.
\end{fact}

The inclusion from right to left follows from the definition of $\LL$. Then we turn to prove the other inclusion. Consider $\langle \A, F\rangle \in \ModstarL$. If $\A$ is trivial, then also $\langle \A, F\rangle = \langle {\bf 1}, \{ 1 \} \rangle$, since $\LL$ has theorems. The suppose that $\A$ is non-trivial. The fact that $\Box \Box x \vdash_{\LL} y$ implies that that characterization of the Leibniz congruence given in point 1 of Lemma \ref{Lem:Polynomial} can be finitized, yielding the following result: for every $p, q \in A$
\begin{equation}\label{Eq:Chap3:Finitized}
p = q \textrm{ if and only if }(p \in F \Leftrightarrow q \in F \textrm{ and }\Box p \in F \Leftrightarrow \Box q \in F).
\end{equation}

Since $\emptyset \vdash_{\LL} 1$ and $\emptyset \vdash_{\LL} \Box 1$, we know that $1, \Box 1 \in F$. The facts that $\Box \Box x \vdash_{\LL} y$ and that $\langle \A, F\rangle$ is non-trivial and reduced, imply that $\Box^{n} 1 \notin F$ for every $n \geq 2$. In particular, this implies that $\Box 1 \ne 1$ and $\Box^{2}1 = \Box^{3}1$ by (\ref{Eq:Chap3:Finitized}). Now let $p \in F$. We have cases: either $\Box p \in F$ or $\Box p \notin F$. By (\ref{Eq:Chap3:Finitized}) in both cases $p \in \{ 1, \Box 1 \}$. Hence $F = \{ 1, \Box 1 \}$.

If $A = \{1, \Box 1, \Box \Box 1\}$, then $\langle \A, F \rangle \cong \langle \A_{3}, \{1, a\}\rangle$. Then consider the case where $A \ne \{ 1, \Box 1, \Box \Box 1 \}$. There is $p \in A \smallsetminus \{ 1, \Box 1, \Box \Box 1 \}$. In particular, this yields that $p \notin F$. By (\ref{Eq:Chap3:Finitized}) and the fact that $p \ne \Box \Box 1$, we know that $\Box p \in F$. Again, since $\Box \Box x \vdash_{\LL} y$ and $\A$ is non-trivial, we obtain that $\Box \Box p \notin F$. Thus from (\ref{Eq:Chap3:Finitized}) we conclude that $\Box p = \Box 1$.

If $A = \{1, \Box 1, \Box \Box 1, p\}$, then $\langle \A, F \rangle \cong \langle \A_{4}, \{ 1, a \} \rangle$. Suppose the contrary towards a contradiction. Then there is $q \in A \smallsetminus \{ 1, \Box 1, \Box \Box 1, c\}$. Since $F = \{ 1, \Box 1 \}$, we know that $q \notin F$. But from (\ref{Eq:Chap3:Finitized}) it follows that either $q = \Box \Box 1$ or $q = p$, against the assumption.
\begin{fact}
The truth sets of $\LL$ are small.
\end{fact}
This is a direct application of the definition of $\LL$ to the characterization of the class $\ModstarL$ given in Fact \ref{Fact:CUIso}.

\begin{fact}
The truth sets of $\LL$ are not equationally definable.
\end{fact}
Observe that the terms in one variable $x$ up to equivalence in $\VVV(\A_{4})$ are $\{ x , \Box x, \Box \Box x, 1, \Box 1 \}$. It is easy to check that $x \thickapprox x$ is the only equation, built up with these terms, that is satisfied by the designated elements of $\langle \A_{4}, \{ 1, a \} \rangle$. Since $\langle \A_{4}, \{ 1, a \} \rangle$ is a reduced model of $\LL$, we conclude that truth is not equationally definable in $\ModstarL$.
\qed
\end{exa}

In Figure \ref{Fig:LeibnizHierarchy2} the reader can find a diagram (where arrows represent inclusions) that subsume the definability conditions considered so far in the framework of the Leibniz hierarchy (see \cite{Cz01,AAL-AIT-f} for the relevant definitions). It is worth to remark that all classes in that diagram are different, as explained in the following technical remark (that the reader may safely skip, without loosing the sense of the paper).

\begin{Remark}
Observe that the fact that the inclusion relations depicted in Figure \ref{Fig:LeibnizHierarchy2} are sound follows from Lemma \ref{Lem:BasicWisdom} and well-known facts about the Leibniz hierarchy. Moreover, to prove that all the classes in the picture are different, it would be enough to show that the $6$ classes below the one of logics, whose truth sets are equationally definable, are different.

The fact that the classes of logics whose truth sets are respectively implicitly definable, small and equationally definable are different is motivated at pag. \pageref{Comments}. Moreover, in Section \ref{Sec:UniExamples} some example of logics whose truth sets are almost parametrically equationally but not equationally, definable are described. Now we describe a logic whose truth sets are almost small, but neither almost parametrically equationally definable nor small. We already motivated the existence of a logic $\LL$ whose truth sets are small, but not equationally definable. Then let $\LL'$ be the logic whose theories are $\Th\LL \cup \{ \emptyset \}$. The truth sets of $\LL'$ are not small, since $\LL$ is purely inferential. Moreover, by Theorem the fact that the truth sets of $\LL$ are not equationally definable,  that $\leibniz^{\A} \colon \FFiLA \to \Con\A$ is not completely-order reflecting for some algebra $\A$. Now, observe that $\FFiLA \subseteq \FFi_{\LL'}A$ and that $\FFiLA$ is a set of non-empty filters. Thus we conclude that the map $\leibniz^{\A} \colon \FFi_{\LL'}\A \to \Con\A$ is not almost completely order-reflecting. By Theorem \ref{Thm:Chap2:LeibnizUniversal} we conclude that the truth sets of $\LL'$ are not almost parametrically equationally definable. Applying a similar argument, it is possible to construct a logic whose truth sets are almost implicitly definable, but neither almost small nor almost implicitly definable. 
\qed
\end{Remark}

\section*{Appendix}

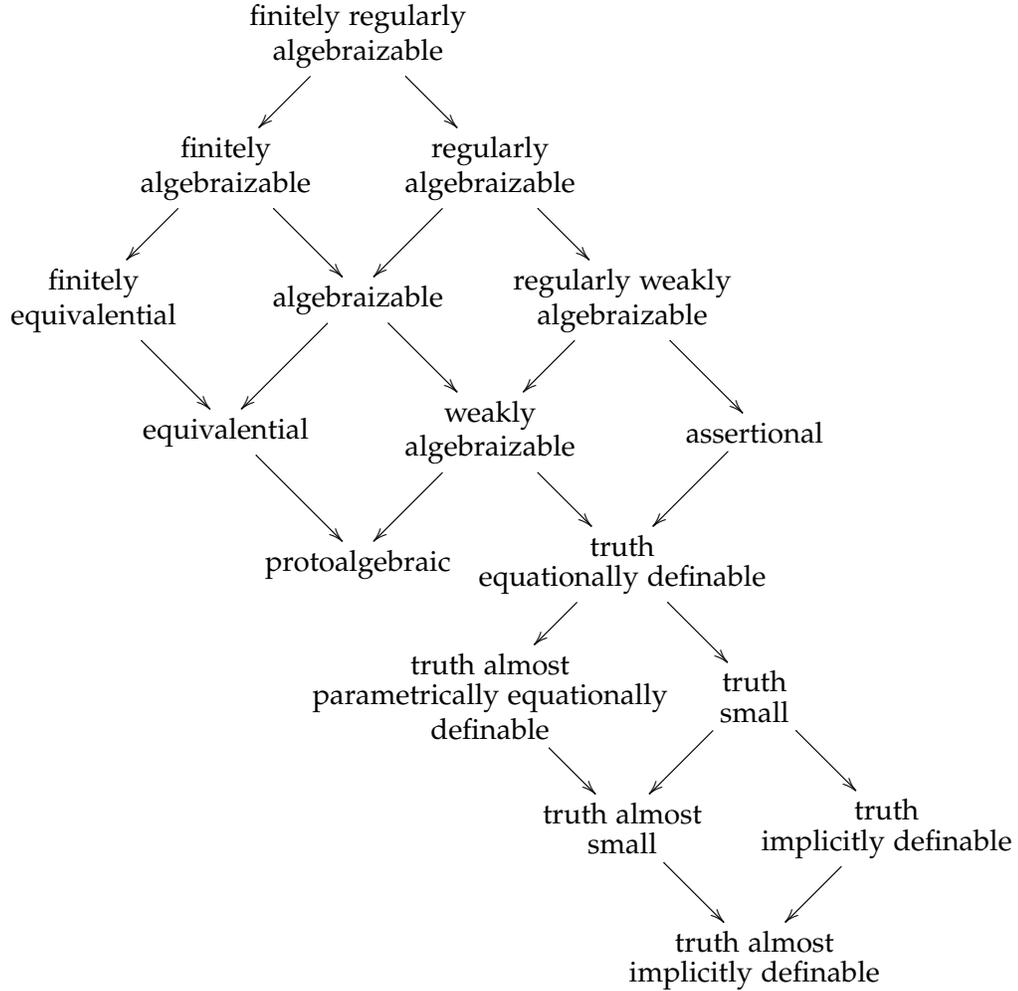
\begin{figure}
\[
\xymatrix@R=50pt @C=50pt @!0{%
 &  & {\txt{finitely regularly \\algebraizable}} \ar[dr]\ar[dl]&&&&\\
 &{\txt{finitely \\algebraizable}} \ar[dr]\ar[dl]  & & {\txt{regularly\\ algebraizable}}\ar[dl]\ar[dr] &&&\\
{\txt{finitely \\equivalential}} \ar[dr]&& {\txt{algebraizable}}\ar[dl]\ar[dr] & & {\txt{regularly weakly\\ algebraizable}}\ar[dl]\ar[dr]&&\\
&{\txt{equivalential}}\ar[dr] & & {\txt{weakly\\ algebraizable}}\ar[dl]\ar[dr] & & {\txt{assertional}}\ar[dl]&\\
& & {\txt{protoalgebraic}} & & {\txt{truth\\ equationally definable}}\ar[dr]\ar[dl]&&\\
& &&{\txt{truth almost\\parametrically equationally\\ definable}}\ar[dr]& &{\txt{truth \\small}}\ar[dl]\ar[dr]&\\
& &&&{\txt{truth almost\\small}} \ar[dr]&&{\txt{truth \\implicitly definable}}\ar[dl]\\
& &&&&{\txt{truth almost\\implicitly definable}}&
}
\]
\caption{The Leibniz hierarchy}\label{Fig:LeibnizHierarchy2}
\end{figure}

In the forthcoming papers \cite{Mor16a,Mor16compl} the problem of classifying logics in the Leibniz hierarchy is studied from a computational point of view. At the time where these papers were written the extension of the Leibniz hierarchy presented here was not developed yet. For this reason, we take the opportunity of mentioning that the proof of \cite[Theorem 4.10]{Mor16a} yields the following:

\begin{Theorem}
Let $\class{K}$ be a level of the Leibniz hierarchy in Figure \ref{Fig:LeibnizHierarchy2}. The problem of determining whether the logic of a given consistent finite Hilbert calculus in a finite language belongs to $\class{K}$ is undecidable.
\end{Theorem}

\

\paragraph{\bfseries Acknowledgements.}

I am very grateful to Josep Maria Font and Ramon Jansana, who read carefully many versions of this work and provided several useful comments that improved its  readability. Thanks are due also to Jos\'{e} Gil-F\'{e}rez and James Raftery for their remarks on a previous manuscript. Finally, I wish to thank the anonymous referees for their helpful comments and suggestions. This research was supported by project GA$17$-$04630$S of the Czech Science Foundation (GA\v{C}R).
\bibliographystyle{plain}

\end{document}